\documentclass[draft,reqno, 12pt]{amsart}
\usepackage{amssymb}
\usepackage{euscript}

\let\cal\mathcal

\makeatletter \@mparswitchfalse \makeatother

\newtheorem{theorem}{Theorem}
\newtheorem{lemma}[theorem]{Lemma}
\newtheorem{corollary}[theorem]{Corollary}
\newtheorem{proposition}[theorem]{Proposition}
\newtheorem{remark}[theorem]{Remark}
\newtheorem{remarks}[theorem]{Remarks}
\newtheorem{definition}[theorem]{Definition}
\newtheorem{varprop}{Proposition}

\numberwithin{equation}{section}
\numberwithin{theorem}{section}

\title[Non-Commutative Martingales]{A Weak type inequality for non-commutative martingales and applications }
\author[N. Randrianantoanina]{Narcisse Randrianantoanina}
\address{Department of Mathematics, Miami
University, Oxford, Ohio 45056  (USA)}
\thanks{Supported in part by NSF grant DMS-0096696.}
\email{randrin@muohio.edu}

\def\N{\mathbb N}

\def\R{\mathbb R}
\def\M{\cal{M}}
\def\H{\cal{H}}
\let\<\langle
\def\ch{\raise 0.5ex \hbox{$\chi$}}
\def\T{\tau}
\def\E{\cal{E}}

\let\\\cr
\let\phi\varphi

\let\epsilon\varepsilon

\def\log{\operatorname{log}}

\begin{document}

\begin{abstract}
We prove a  weak-type (1,1) inequality for  square functions of
 non-commutative  martingales that are simultaneously bounded in $L^2$ and $L^1$.
 More precisely, the
 following non-commutative analogue of a classical result of
 Burkholder holds:
   there exists an absolute  constant $K>0$ such that
if $\M$ is a  semi-finite von Neumann algebra and
$(\M_n)^{\infty}_{n=1}$ is an increasing filtration of von Neumann
subalgebras of $\M$ then  for any given  martingale
$x=(x_n)^{\infty}_{n=1}$ that is    bounded in $L^2(\M)\cap
L^1(\M)$, adapted to $(\M_n)^{\infty}_{n=1}$,  there exist two
\underline{martingale difference} sequences,
$a=(a_n)_{n=1}^\infty$ and $b=(b_n)_{n=1}^\infty$, with $dx_n =
a_n + b_n$ for every $n\geq 1$,
\[
\left\| \left(\sum^\infty_{n=1} a_n^*a_n
\right)^{{1}/{2}}\right\|_{2} + \left\| \left(\sum^\infty_{n=1}
b_nb_n^*\right)^{1/2}\right\|_{2} \leq  2\left\| x \right\|_2,
\]
and
\[
\left\| \left(\sum^\infty_{n=1} a_n^*a_n
\right)^{{1}/{2}}\right\|_{1,\infty} + \left\|
\left(\sum^\infty_{n=1} b_nb_n^*\right)^{1/2}\right\|_{1,\infty}
\leq K\left\| x \right\|_1.
\]

 As an application, we obtain the optimal orders  of growth for the
constants involved in the Pisier-Xu  non-commutative analogue of
the classical Burkholder-Gundy inequalities.
\end{abstract}

\maketitle

\setcounter{section}{-1}

\section{Introduction}

Non-commutative (or quantum) probability has developed
considerably in recent years. It provides many connections between
several fields of mathematics such as mathematical physics,
operator algebras, and classical probability theory. We refer to
the book by Meyer \cite{Mey} for general quantum probability, the
book by Parthasarathy \cite{Part}  for quantum stochastic
calculus, and the book by Voiculescu et al. \cite{VDN} for free
probability.

In classical theory, martingale theory has played a significant
role in the developments of various  fields of analysis (see for
instance \cite{Bu3,EG,LT}). In this paper, our main interest is on
non-commutative martingales. As in the classical case,
non-commutative martingales have connections with other area such
as operator algebra theory, operator space theory, and matrix
valued harmonic analysis which includes  among other things,
operator valued Carleson measures, operator valued Hardy spaces,
and operator valued Hankel operators (see for instance
\cite{GPTV,NPTV}).

 Alongside  the general development  of quantum probability
theory, the subfield of non-commutative martingales has received
considerable progress in recent years. Indeed, many of  classical
inequalities from the usual (commutative) martingale theory have
been generalized to the non commutative settings. Let us recall
some sample contributions by several authors.  For instance,
pointwise convergence of non-commutative martingales was already
considered by Dang-Ngok \cite{DN}, Cuculescu \cite{Cuc}, and
Barnett \cite{BA} in the 70's and 80's. Pisier and Xu \cite{PX}
proved the non-commutative analogue of the Burkholder-Gundy
inequalities on square functions and non commutative analogue of
Stein's inequality. It is their general functional analytic
approach that led to the consideration of non-commutative analogue
of several classical martingale inequalities.  A non-commutative
analogue of Doob's maximal inequality was successfully  formulated
and proved by Junge in \cite{Ju} and non-commutative analogues of
Burkholder/Rosenthal inequalities on conditioned square functions
were studied  by Junge and Xu in \cite{JX} among many other
related topics. These different results pave the way to the
consideration of non-commutative martingale Hardy spaces and
non-commutative martingale $BMO$ which are non-commutative
generalizations of spaces  that are central to the developments of
classical harmonic analysis and interpolation theory. We note also
a very recent result of Musat \cite{MUS} on interpolation
involving non-commutative $BMO$ and non-commutative $L^p$-spaces
as endpoints.

In most of the papers listed above, square functions played a very
crucial role. Note however that in strong
 contrast with the classical case, square functions in the
non-commutative case can take many different forms so it is very
important to formulate the
 ``right" square functions. Recall that if
$1<p<2$, and $x=(x_n)_{n=1}^\infty$ is a non-commutative
martingale (see the formal definition below), the $\cal{H}^p$-norm
($\cal{H}^p$ being the Hardy space of non-commutative martingales)
 is given by:
\begin{equation}\label{square}
\|x\|_{\cal{H}^p}=\inf\left\{\left\|\left(\sum_{n\geq 1}
|dy_n|^2\right)^{{1}/{2}}\right\|_p + \left\|\left(\sum_{n\geq 1}
|dz_n^*|^2\right)^{{1}/{2}}\right\|_p\right\}
\end{equation}
where the infimum runs over all decompositions $x=y+z$, with $y$
and $z$ being martingales. That is, it depends on two different
types of square functions (given by right and left moduli). The
fact that one has to decompose the martingale $x$ into two
martingales was first discovered for non-commutative Khintchine
inequalities (see \cite{LP4,LPI}) and this type of decomposition
is often the source of the difficulties in extending classical
results to non-commutative settings.

 The main purpose of this paper is  to study the square functions
 of non commutative martingales for
 the case $p=1$ which is primarily motivated by the following classical
 result of Burkholder.
\begin{theorem}[\cite{Bu4}]
 Let $(f_n)^\infty_{n=1}$ be a martingale on a
probability space $(\Omega,\Sigma, P)$ and
$S(f)=(\sum_{n=1}^\infty |f_n -f_{n-1}|^2)^{1/2}$. Then there
exists an absolute constant $M>0$ such that for every $\lambda
>0$,
\[
\lambda P\left( S(f)>\lambda\right) \leq M\sup_n \mathbb{E
}(|f_n|).
\]
\end{theorem}
It is a natural question to consider whether  Theorem~0.1 has
non-commutative counterparts. We remark that Burkholder deduced
the above result from  the weak-type (1,1) boundedness of
martingale transforms (also proved in \cite{Bu4}) via the
classical Khintchine inequality. One can also use  the  classical
Doob's identity  (see for instance \cite[Chap.~II]{GA}) to deduce
Theorem~0.1 from the weak-type (1,1) boundedness of  martingale
transforms.
 We note that non-commutative martingale transforms
are of weak-type (1,1) (\cite{Ran15}). However, unlike the
classical case, a non-commutative analogue of Theorem 0.1 can not
be deduced directly from the weak type (1,1) boundedness of
martingale transforms via the classical techniques, as (at least
at the time of this writing) there is no adequate Khintchine
inequality for non-commutative weak-$L^1$-spaces.

In \cite{Ran16}, a first attempt was made to generalize
Theorem~0.1 to  non-commutative settings. The $p$-norm of the
square functions for the case $1<p<2$  stated in (\ref{square})
suggests that  the formulation of the weak-$L^1$ norm  of square
functions should require decompositions of the martingales
involved. We obtained in \cite[Theorem~2.1]{Ran16}  a
decomposition of any given martingale into two sequences in weak
$L^1$-space where the corresponding weak $L^1$-norm of the square
functions similar to the one stated in $(\ref{square})$ is bounded
by the $L^1$-norm of the corresponding martingale. The result from
\cite{Ran16} prompted the question of whether or not such
decomposition can be chosen to be martingales. For the finite
case, our main result answers this positively for the case of
$L^2$-bounded martingales (see Theorem~\ref{main} below). More
precisely, there exists a constant $K>0$ such that if
$x=(x_n)_{n=1}^\infty$ is a non-commutative martingale, there
exists two martingales $y=(y_n)_{n=1}^\infty$ and
$z=(x_n)_{n=1}^\infty$ such that $x=y+z$ and with the property:
\begin{equation}
\left\|\left(\sum_{n\geq 1}
|dy_n|^2\right)^{{1}/{2}}\right\|_{1,\infty} +
\left\|\left(\sum_{n\geq 1}
|dz_n^*|^2\right)^{{1}/{2}}\right\|_{1,\infty} \leq
\left\|x\right\|_1.
\end{equation}
 Coupled
with general interpolation techniques, our main result provides a
solution to a problem left open in \cite{JX} (see also
\cite[Problem~8.2]{X4}) on optimal order of growth of the
constants involved in the non-commutative Burkholder-Gundy
inequalities when $p \to 1$ (see Theorem~\ref{alpha} below).

In order to achieve the decomposition into two martingales, our
method of proof (although it follows closely those taken in
\cite{Ran16} and \cite{Ran15})
 requires substantial adjustments.
 It  depends  heavily on a non-commutative version
of the classical Doob's maximal inequality obtained by Cuculescu
\cite{Cuc},  and weak-type (1,1) boundedness of triangular
truncations relative to disjoint projections.

The paper is organized as follows:  in Section~1 below, we set
some basic preliminary background concerning  non-commutative
spaces  and collect some results on triangular truncations. In
Section~2, we recall the general setup of  non-commutative
martingale theory. Section~3 is devoted mainly to the statement of
the main decomposition, the construction of the decomposition and
a detailed proof of the weak-type inequality  for the finite case.
In Section~4,  we  will point out the adjustment needed to extend
our main result  from Section~3 to the semi-finite case. In the
last section, we provide a new proof of one of the inequalities
involved in the non-commutative Burkholder-Gundy inequality and
deduce the optimal order of the constants involved.

\section{Non-commutative spaces and preliminary results }
We use standard notation in operator algebras. We refer to
\cite{KR} and \cite{TAK} for background  on von Neumann algebra
theory. In this section, we will recall some basic definitions
that we will use throughout this paper. In particular, we will
outline the general construction of non-commutative spaces and
discuss triangular truncations with respect to sequence of
disjoint projections.

 Throughout,  $\M$ is  a semi-finite von Neumann
algebra with a normal faithful semi-finite  trace  $\T$. The
identity element of $\M$ is denoted by ${\bf 1}$. For $0 < p\leq
\infty$, let $L^p(\M,\T)$ be the associated non-commutative
$L^p$-space (see for instance \cite{DIX} and \cite{N}). Note that
if $p=\infty$, $L^\infty(\M,\tau)$ is just $\M$ with the usual
operator norm; also recall that for $0< p<\infty$, the
(quasi)-norm on $L^p(\M,\T)$ is defined by
$$\Vert x \Vert_p =( \T(|x|^p))^{1/p}, \qquad x \in L^p(\M,\T),$$
where $|x|=(x^*x)^{1/2}$ is the usual modulus of $x$.

In order to ease the introduction of some of the spaces used in
the sequel, we need the general scheme of symmetric spaces of
measurable operators developed in \cite{CS,DDP1,DDP3,X}.

Let $H$ be a  complex Hilbert space and $\M \subseteq B(H)$.
 A closed densely defined operator $a$ on $H$ is said to be {\em
affiliated with} $\M$ if $u^* au = a$ for all unitary $u$ in the
commutant $\M'$ of $\M$. If $a$ is a densely defined self-adjoint
operator on $H$, and if $a = \int^\infty_{- \infty} s d e^a_s$ is
its spectral decomposition, then for any Borel subset $B \subseteq
\R$, we denote by $\ch_B(a)$ the corresponding spectral projection
$\int^\infty_{- \infty} \ch_B(s) d e^a_s$. A closed densely
defined operator $a$ on $H$ affiliated with $\M$ is said to be
{\em $\T$-measurable} if there exists a number $s \geq 0$ such
that $\T(\ch_{(s, \infty)} (|a|)) < \infty$.

The set of all $\T$-measurable operators will be denoted by
$\overline{\M}$. The set $\overline{\M}$ is a $*$\ -algebra with
respect to the strong sum, the strong product, and the adjoint
operation \cite{N}.  For $\epsilon, \delta>0$, let
\[
N(\epsilon,\delta)=\{x\in \overline{\M} : \mathrm{for\ some\
projection}\ p\in \M, \Vert xp\Vert<\epsilon\ \mathrm{and}\
\T({\bf 1}-p)\leq \delta\}.
\]
The system $(N,\epsilon,\delta))_{\epsilon,\delta}$ forms a
fundamental system of neighborhoods of the origin of the vector
space $\overline{\M}$ and the translation-invariant topology
induced by this system is called the {\it measure topology}.
Convergence in measure will be used in the sequel.

For $x \in \overline{\M}$, the generalized singular value function
$\mu (x)$ of $x$ is defined by
\[
\mu_t(x) = \inf \{ s \geq 0: \T(\ch_{(s, \infty)} (|x|)) \leq t
\}, \quad \text{ for }\  t \geq 0.
\]
The function $t \to \mu_t(x)$ from  the interval $[0, \T({\bf
1}))$ to $[0, \infty)$ is right continuous, non-increasing and is
the inverse of the distribution function $\lambda (x)$, where
$\lambda_s(x) = \T(\ch_{(s, \infty)}(|x|))$, for $s \geq 0$. For
an in depth study  of $\mu(.)$ and $\lambda(.)$,  we refer the
reader to \cite{FK}.

 For the definition below, we refer
the reader to \cite{BENSHA,LT} for the theory of rearrangement
invariant function spaces.

\begin{definition}
Let $E$ be a rearrangement invariant (quasi-) Banach function
space on  the interval $[0, \T({\bf 1}))$. We define the symmetric
space $E(\M, \T)$ of measurable operators by setting:
\begin{eqnarray*}
E(\M, \T) &=& \{ x \in
\overline{\M}\ : \ \mu(x) \in E \} \quad \text{and} \\
\|x\|_{E(\M,\T)} &=& \| \mu(x)\|_E,\  \text{for}\  x \in E(\M,\T).
\end{eqnarray*}
\end{definition}

It is well known that $E(\M, \T)$ is a Banach space (respectively,
quasi-Banach space) if $E$ is a Banach space (respectively,
quasi-Banach space). The space $E(\M, \T)$ is often referred to as
the non-commutative analogue of the function space $E$ and
  if $E = L^p[0,\T({\bf 1}))$, for $0 < p \leq \infty$, then
$E(\M, \T)$ coincides with the usual non-commutative $L^p$\ -space
associated with $(\M, \T)$. We refer to \cite{CS,DDP1,DDP3,X}
 for more detailed discussions about these
spaces. Of special interest in this paper are non-commutative weak
$L^1$-spaces. The non-commutative weak $L^1$-space, denoted by
$L^{1,\infty}(\M,\T)$,  is defined as the linear subspace of all
$x \in \overline{\M}$ for which  the quasi-norm
\begin{equation}\label{weaknorm}
\Vert x \Vert_{1,\infty} := \sup_{t>0} t\mu_t(x) =
\sup_{\lambda>0}\lambda \T(\ch_{(\lambda,\infty)}(|x|))
\end{equation}
is finite. Equipped with the quasi-norm $\Vert \cdot
\Vert_{1,\infty}$, $L^{1,\infty}(\M,\T)$ is a quasi-Banach space.
It is easy to verify that as in the commutative space,  $  \Vert x
\Vert_{1,\infty} \leq \Vert x\Vert_1 $ for all $x \in L^1(\M,\T)$.

For a complete, detailed, and  up to date presentation of
non-commutative integration and non-commutative spaces, we refer
to the recent survey article \cite{PX3}.

The next lemma  is probably well known. It will be used repeatedly
in the sequel.
\begin{lemma}\label{splitting}
Let $a$ and $b$ be operators in $L^{1,\infty}(\M,\T)$. For every
$\lambda>0$, $\alpha\in (0,1)$, and $\beta \in (0,1)$,
\[\T(\ch_{(\lambda, \infty)}(|a+b|))\leq \alpha^{-1}
\T(\ch_{(\beta\lambda, \infty)}(|a|)) + (1-\alpha)^{-1}
\T(\ch_{((1-\beta)\lambda, \infty)}(|b|)).
\]
\end{lemma}
\begin{proof}
  Using properties of
generalized singular value functions $\mu (\cdot)$ from \cite{FK},
we have,
\[ \T\left(\ch_{(\lambda,\infty)}(|a+b|)\right) =
\int^1_0 \ch_{(\lambda, \infty)} \{\mu_t(a+b)\}\ dt.
\]
This follows from \cite[Corollary~2.8]{FK}  by approximating  the
characteristic function $\ch_{(\lambda, \infty)}(\cdot)$ from
below by sequences of continuous functions $f$ on $[0,\infty)$
satisfying $f(0)=0$. We can deduce the following estimate:
\begin{equation*}
\begin{split}
\T\left(\ch_{(\lambda,\infty)} (|a+b|)\right) &\leq \int^1_0
\ch_{(\lambda,\infty)} \{\mu_{\alpha{t}}(a)  +
\mu_{(1-\alpha)t}(b)\}\
dt\\
&\leq \int^1_0 \ch_{({\beta}\lambda,\infty)}
\{\mu_{\alpha{t}}(a)\}\ dt + \int^1_0
\ch_{((1-{\beta})\lambda,\infty)} \{\mu_{(1-\alpha){t}}(b)\}\ dt
\end{split}
\end{equation*}
and by  simple change of variables,
\begin{equation*}
\begin{split}
\T\left(\ch_{(\lambda,\infty)}(|a+b|)\right) &\leq \alpha^{-1}
\int^1_0 \ch_{({\beta}\lambda,\infty)} \{\mu_{t}(a)\}dt +
(1-\alpha)^{-1}
\int^1_0 \ch_{((1-{\beta})\lambda,\infty)}\{\mu_{t}(b)\} dt \\
&=\alpha^{-1} \T(\ch_{(\beta\lambda, \infty)}(|a|)) +
(1-\alpha)^{-1} \T(\ch_{((1-\beta)\lambda, \infty)}(|b|))
\end{split}
\end{equation*}
as stated in the lemma.
\end{proof}


We end this section with a brief  discussion on {\it triangular
truncations}. This will be very crucial throughout the paper. Let
$\cal{P}=\{p_i\}_{i=1}^M$ be an arbitrary finite or infinite
sequence of mutually orthogonal projections from $\M$. We recall
the triangular truncation on $\overline{\M}$ (with respect to
$\cal{P}$) by
\[
\cal{T}^{(\cal{P})}x :=\sum_{j=1}^M\sum_{i\leq j} p_ixp_j, \quad x
\in \overline{\M}.
\]
The diagonal projection $D^{(\cal{P})}$ is defined on
$\overline{\M}$ by setting
\[
D^{(\cal{P})}x := \sum_{i=1}^M  p_ixp_i, \quad x \in
\overline{\M}.
\]
We also use the following operator on $\overline{\M}$,
\[
\cal{H}^{(\cal{P})}x:= -i(\cal{T}^{(\cal{P})}x -
\cal{T}^{(\cal{P})}x^*), \quad x \in \overline{\M}.
\]
For convenience, we collect some properties of the operators
introduced above that are useful for our presentation. In the
following lemma, $\overline{\M}_{\cal{P}}$ denotes the range of
$\cal{T}^{(\cal{P})}$.
\begin{lemma}[\cite{DDPS2}]\label{diagonal}
The  operators defined above satisfy the following properties:
\begin{itemize}
\item[(i)] If $0\leq x \in \overline{\M}$ and $\lambda>0$, then
$\lambda{\bf 1} +x + i\cal{H}^{(\cal{P})}(x)$ is invertible, with
$(\lambda{\bf 1} +x + i\cal{H}^{(\cal{P})}(x))^{-1} \in \M$ and
$\|(\lambda{\bf 1} +x + i\cal{H}^{(\cal{P})}(x))^{-1}\|_{\infty}
\leq 1/\lambda$.
\item[(ii)]
$D^{(\cal{P})}\cal{T}^{(\cal{P})}=\cal{T}^{(\cal{P})}D^{(\cal{P})}=D^{(\cal{P})}$.
\item[(iii)] If $x$, $y\in \overline{\M}_{\cal{P}}$, then
$D^{(\cal{P})}(xy)= D^{(\cal{P})}(x)D^{(\cal{P})}(y)$.
\end{itemize}
If we assume that $\sum_{i=1}^M  p_i={\bf 1}$, then:
\begin{itemize}
\item[(iv)] If $x \in \overline{\M}_{\cal{P}}$ is invertible and
$ D^{(\cal{P})}(x)$ is self-adjoint, then $x^{-1} \in
\overline{\M}_{\cal{P}}$  and $D^{(\cal{P})}(x)^{-1}=
D^{(\cal{P})}(x^{-1})$.
\item[(v)] If $x \in \overline{\M}$ is self-adjoint, then
$x +i\cal{H}^{(\cal{P})}(x)=2\cal{T}^{(\cal{P})}(x)$.
\item[(vi)] If $x \in L^1(\M,\T)$, then
$\T(D^{(\cal{P})}(x))=\T(x)$.
\end{itemize}
\end{lemma}

The next result is a weak-type boundedness  of \lq\lq$l^2$-sum" of
finite family of triangular truncations.

\begin{proposition}\label{truncation}
Let $\{\cal{P}^{(n)}\}_{n=1}^N$ be a family of finite sequence of
mutually disjoint projections with
$\cal{P}^{(n)}=\{p_{i,n}\}_{i=1}^M$ for each $1\leq n\leq N$. If
$(x_n)_{n=1}^N$ is a finite sequence of positive operators in
$L^1(\M,\T)$ then
\[
\left\| \left(\sum_{n=1}^N
|\cal{T}^{(\cal{P}^{(n)})}x_n|^2\right)^{1/2}\right\|_{1,\infty}
\leq  5\sqrt{2} \sum_{n=1}^N \|x_n\|_1.
\]
\end{proposition}

We remark that if $N=1$, then the above proposition (at least for
the finite case)  is a particular case of the weak type (1,1)
boundedness of the Hilbert transform associated with finite
subdiagonal subalgebra obtained in \cite{Ran9}. A  more concise
proof for the case $N=1$ also appeared in
\cite[Theorem~1.4]{DDPS2}. It is the presentation in \cite{DDPS2}
that we will  adopt below to prove Proposition~\ref{truncation}.


\begin{proof}[Proof of Proposition~\ref{truncation}] For each $1\leq n\leq N$,
we will simply write $\cal{T}_n$, $D_n$, and $\cal{H}_n$ for
$\cal{T}^{(\cal{P}^{(n)})}$, $D^{(\cal{P}^{(n)})}$, and
$\cal{H}^{(\cal{P}^{(n)})}$, respectively.

Since $\|(\sum_{i=1}^M p_{i,n})x(\sum_{i=1}^M p_{i,n})\|_1\leq
\|x\|_1$ and $\cal{T}_n((\sum_{i=1}^M p_{i,n})x(\sum_{i=1}^M
p_{i,n}))=\cal{T}_nx$ for every $x \in L^1(\M,\T)$, it is clear
that it is enough to consider the case where  for each $1\leq
n\leq N$, $x_n $ belongs to the space $ L^1((\sum_{i=1}^M
p_{i,n})\M(\sum_{i=1}^M p_{i,n}))$. Therefore we may assume
without loss of generality that  for each $1\leq n\leq N$,
$\sum_{i=1}^M p_{i,n}={\bf 1}$.  We will assume first that for
each $1\leq n\leq N$, $x_n \in \M \cap L^1(\M,\T)$.

For $1\leq n\leq N$, set
\[
A_n:= x_n +i \cal{H}_n(x_n).
\]
We will show first that for $\lambda>0$,
\begin{equation}\label{hilbert}
\T\left(\ch_{(\lambda, \infty)}((\sum_{n=1}^N
|A_n|^2)^{1/2})\right)\leq 4 (\sum_{n=1}^N \|x_n\|_1)/\lambda.
\end{equation}

 The main argument is to estimate the trace of the operator
$\sum_{n=1}^N |D_n(A_n(\lambda{\bf 1}+A_n)^{-1})|$ for $\lambda>0$
from above and below.

Note that for $1\leq n\leq N$, $A_n \in
\overline{\M}_{\cal{P}^{(n)}}$ and $D_n(A_n)=D_n(x_n)$. We can
deduce from Lemma~\ref{diagonal}(iii) that
\[D_n(A_n(\lambda{\bf
1}+A_n)^{-1})=D_n(A_n)D_n((\lambda{\bf 1}+A_n)^{-1}).\] The
following estimate from above   follows  directly from
Lemma~\ref{diagonal}(i):
\begin{equation}\label{above}
\begin{split}
\sum_{n=1}^N \T(|D_n(A_n(\lambda{\bf 1}+A_n)^{-1})|)&\leq
\sum_{n=1}^N \T(D_n(x_n))\|(\lambda{\bf 1}+A_n)^{-1}\|_\infty \\
&= \sum_{n=1}^N \|x_n\|_1 . \|(\lambda{\bf 1}+A_n)^{-1}\|_\infty \\
&\leq (\sum_{n=1}^N \|x_n\|_1 )/\lambda.
\end{split}
\end{equation}

For the estimate from below, we first note from
Lemma~\ref{diagonal}(iii), (iv), and (v) that for $1\leq n\leq N$,
the operator
\[
D_n(A_n(\lambda{\bf 1}+A_n)^{-1})=D_n(A_n)D_n((\lambda{\bf
1}+A_n)^{-1})= D_n(x_n)(\lambda{\bf 1}+D_n(x_n))^{-1}
\]
is self-adjoint and we clearly have,
\[
\sum_{n=1}^N \T(|D_n(A_n(\lambda{\bf 1}+A_n)^{-1})|)\geq
\sum_{n=1}^N \T(D_n(A_n(\lambda{\bf 1}+A_n)^{-1})).
\]
For each $1\leq n\leq N$, we will estimate $\T(D_n(A_n(\lambda{\bf
1}+A_n)^{-1}))$ exactly as in \cite{DDPS2}. We include the
argument for completeness.
\begin{eqnarray*}
\T(D_n(A_n(\lambda{\bf
1}+A_n)^{-1}))&=&\T(\mathrm{Re}D_n(A_n(\lambda{\bf
1}+A_n)^{-1}))\\
&=&\T(D_n(\mathrm{Re}A_n(\lambda{\bf
1}+A_n)^{-1}))\\
&=&\T(\mathrm{Re}A_n(\lambda{\bf 1}+A_n)^{-1}).
\end{eqnarray*}
Note  that $\mathrm{Re}A_n(\lambda{\bf 1}+A_n)^{-1}=(\lambda{\bf
1}+A_n^*)^{-1}(|A_n|^2 +\lambda\mathrm{Re}A_n)(\lambda{\bf 1}
+A_n)^{-1}\geq 2^{-1}[(\lambda{\bf 1}+A_n^*)^{-1}(|A_n|^2
+2\lambda\mathrm{Re}A_n)(\lambda{\bf 1} +A_n)^{-1}]$. We have,
\begin{eqnarray*}
\T(D_n(A_n(\lambda{\bf 1}+A_n)^{-1})) &\geq&
2^{-1}\T\left((\lambda{\bf 1}+A_n^*)^{-1}(|A_n|^2
+2\lambda\mathrm{Re}A_n)(\lambda{\bf 1} +A_n)^{-1}\right)\\
&=& 2^{-1}\T\left((|A_n|^2
+2\lambda\mathrm{Re}A_n)(\lambda{\bf 1} +A_n)^{-1}(\lambda{\bf 1}+A_n^*)^{-1}\right)\\
&=& 2^{-1}\T\left((|A_n|^2 +2\lambda\mathrm{Re}A_n)(|A_n|^2
+2\lambda\mathrm{Re}A_n +\lambda^2{\bf 1})^{-1}\right).
\end{eqnarray*}
Set $y_n:=|A_n|^2 +2\lambda\mathrm{Re}A_n$ and $y:=\sum_{n=1}^N
y_n$. As $\mathrm{Re}A_n\geq 0$, we have $y_n\geq 0$ and therefore
for each $n\geq 1$,  $y_n\leq y$ and $(y_n +\lambda^2{\bf
1})^{-1}\geq (y +\lambda^2{\bf 1})^{-1}$. Hence, we can deduce,
\begin{equation*}
\begin{split}
\sum_{n=1}^N \T(D_n(A_n(\lambda{\bf 1}+A_n)^{-1})) &\geq
2^{-1}\sum_{n=1}^N \T(y_n(y_n +\lambda^2{\bf 1})^{-1})\\
&\geq 2^{-1}\sum_{n=1}^N \T(y_n(y +\lambda^2{\bf 1})^{-1})\\
&=2^{-1} \T(y(y +\lambda^2{\bf 1})^{-1}).
\end{split}
\end{equation*}
Observe that $y(y +\lambda^2{\bf 1})^{-1}\geq \ch_{(\lambda^2,
\infty)}(y). y(y +\lambda^2{\bf 1})^{-1}\geq
2^{-1}\ch_{(\lambda^2, \infty)}(y)$. We obtain,
\[ \sum_{n=1}^N
\T(D_n(A_n(\lambda{\bf 1}+A_n)^{-1})) \geq 4^{-1}
\T(\ch_{(\lambda^2, \infty)}(y)).
\]
Again as $\mathrm{Re}A_n\geq 0$, $y \geq \sum_{n=1}^N |A_n|^2$, we
have
\begin{equation}\label{below}
\begin{split}
\sum_{n=1}^N \T\left(D_n(A_n(\lambda{\bf 1}+A_n)^{-1})\right)
&\geq 4^{-1} \T\left(\ch_{(\lambda^2, \infty)}(\sum_{n=1}^N
|A_n|^2)\right)\\
&=  4^{-1} \T\left(\ch_{(\lambda, \infty)}((\sum_{n=1}^N
|A_n|^2)^{1/2})\right).
\end{split}
\end{equation}
Combining (\ref{above}) and (\ref{below}),
inequality~(\ref{hilbert}) follows and hence
\[
\left\|(\sum_{n=1}^N |A_n|^2)^{1/2}\right\|_{1,\infty}\leq 4
\sum_{n=1}^N \left\|x_n\right\|_1.
\]

Now, from Lemma~\ref{diagonal}(vi) and the  elementary fact that
for any operators $a$ and $b$, $|a+b|^2\leq 2(|a|^2 +|b|^2)$, we
have:
\[
\sum_{n=1}^N |\cal{T}_n(x_n)|^2 \leq 2^{-1}(\sum_{n=1}^N |A_n|^2 +
\sum_{n=1}^N |D_n(x_n)|^2).
\]
Using properties of singular values $\mu(\cdot)$ from \cite{FK},
we have   for $t>0$,
\begin{equation*}
\begin{split} t\mu_t\left\{(\sum_{n=1}^N
|\cal{T}_n(x_n)|^2)^{1/2}\right\}&\leq
\sqrt{2}\frac{t}{2}\mu_{t/2}\left\{(\sum_{n=1}^N
|A_n|^2)^{1/2}\right\}\\
&\  + \sqrt{2}\frac{t}{2}
\mu_{t/2}\left\{(\sum_{n=1}^N |D_n(x_n)|^2)^{1/2}\right\} \\
&\leq \sqrt{2}\|(\sum_{n=1}^N |A_n|^2)^{1/2}\|_{1,\infty} +
\sqrt{2}\|(\sum_{n=1}^N |D_n(x_n)|^2)^{1/2}\|_{1,\infty}\\
&\leq 4\sqrt{2}\sum_{n=1}^N\|x_n\|_1 + \sqrt{2}\|(\sum_{n=1}^N
|D_n(x_n)|^2)^{1/2}\|_{1}.
\end{split}
\end{equation*}
Note  that $\|(\sum_{n=1}^N |D_n(x_n)|^2)^{1/2}\|_{1}\leq
\sum_{n=1}^N \|D_n(x_n)\|_1 \leq  \sum_{n=1}^N \|x_n\|_1$. It
follows that for every $t>0$,
\[
t\mu_t\left\{(\sum_{n=1}^N |\cal{T}_n(x_n)|^2)^{1/2}\right\} \leq
5\sqrt{2}\sum_{n=1}^N \|x_n\|_1.
\]
 Taking the supremum over $t>0$, we have proved that the
proposition holds for all finite sequences of positive operators
in $\M \cap L^1(\M,\T)$.

 We complete the proof of the proposition  by noting that if
$(x_n)_{n=1}^N$ is a finite sequence of positive operators in
$L^1(\M,\T)$ then for each $1\leq n\leq N$, we can choose a
sequence $(x_{n}^{(k)})_{k=1}^\infty$ in $\M \cap L^1(\M,\T)$ with
$0\leq x_{n}^{(k)} \uparrow^k x_n$. Observe that for every $1\leq
n\leq N$,  $\cal{T}_n(x_{n}^{(k)}) \rightarrow_k \cal{T}_n(x_n)$
in $L^{1,\infty}(\M,\T)$.  A fortiori, the sequence
$\{(\sum_{n=1}^N|\cal{T}_n(x_{n}^{(k)})|^2)^{1/2}\}_{k=1}^\infty$
converges to   $(\sum_{n=1}^N|\cal{T}_n(x_{n})|^2)^{1/2}$ for the
measure topology (when $k \to \infty$). From \cite[Lemma~3.1]{FK},
we can conclude that   for every $t>0$,
 \[\lim_{k\to
\infty}\mu_t\{(\sum_{n=1}^N|\cal{T}_n(x_{n}^{(k)})|^2)^{1/2}\}=
\mu_t\{(\sum_{n=1}^N|\cal{T}_n(x_{n})|^2)^{1/2}\}.\]
 Hence, for
$t>0$,
\begin{equation*}
\begin{split}
t\mu_t\left\{(\sum_{n=1}^N|\cal{T}_n(x_{n})|^2)^{1/2}\right\}&=
\lim_{k\to\infty}t\mu_t\left\{(\sum_{n=1}^N|\cal{T}_n(x_{n}^{(k)})|^2)^{1/2}\right\}\\
&\leq 5\sqrt{2} \lim_{k\to \infty} \sum_{n=1}^N \|x_{n}^{(k)}\|_1 \\
&= 5\sqrt{2}\sum_{n=1}^N \|x_{n}\|_1.
\end{split}
\end{equation*}
Taking the supremum over $t>0$, the definition of
$\|\cdot\|_{1,\infty}$ provides the desired inequality and thus
the proof of the proposition is complete.
\end{proof}

\section{Conditional expectations and non-commutative martingales}

Let $(\M,\T)$ be  a semi-finite von Neumann algebra and $\cal{N}$
be a von Neumann subalgebra of $\M$. A  linear map $\E: \M \to
\cal{N}$ is  called a {\em normal conditional expectation} if it
satisfies the following:
\begin{itemize}
\item[(i)] $\E$ is a weak$^*$-continuous  projection;
\item[(ii)] $\E$ is positive;
\item[(iii)] $\E(axb)=a\E(x)b$ for all $a, b \in \cal{N}$ and $x
\in \M$;
\item[(iv)] $\T\circ\E=\T$.
\end{itemize}

Recall that such normal conditional expectation from $\M$ onto
$\cal{N}$  exists if and only if the restriction of the trace of
$\M$ to $\cal{N}$ remains semi-finite. For the case where $\M$ is
finite, such conditional expectations always exist. Indeed, if
  $\cal{N}$ is a von Neumann subalgebra of a finite von Neumann algebra $\M$,
   then  the
embedding $\iota: L^1(\cal{N},\T|_{\cal{N}}) \to L^1(\M,\T)$ is an
isometry and the dual map $\E=\iota^*: \M \to \cal{N}$ yields a
 normal conditional expectation (see for instance,
\cite[Theorem~3.4]{TAK}).

Since $\E$ is trace preserving,  it extends as a contractive
projection $\E: L^p(\M,\T)\to L^p(\cal{N},\T|_{\cal{N}})$ for all
$1\leq p \leq \infty$   satisfying the property:
$$ \E(axb)=a\E(x)b$$
when $1\leq p,\ q,\ r\leq \infty$, $1/p +1/q +1/r \leq 1$, $a\in
L^p(\cal{N},\T|_{\cal{N}})$, $b\in L^q(\cal{N},\T|_{\cal{N}})$ and
$x \in L^r(\M,\T)$. More generally, a simple interpolation
argument would prove that if $E$ is a rearrangement invariant
Banach function space on $[0, \T({\bf 1}))$, then $\E$ is a
contraction from $E(\M,\T)$ onto $E(\cal{N},\T|_{\cal{N}})$.

Let us  recall the general setup for martingales. The reader is
referred to \cite{Doob} and \cite{GA} for the classical
(commutative) martingale theory. Let $(\M_n)_{n=1}^\infty$ be an
increasing sequence of von Neumann subalgebras of $\M$ such that
the union of $\M_n$'s is weak$^*$-dense in $\M$. For each $n\geq
1$, assume that there  is a normal conditional expectation $\E_n$
 from $\M$ onto $\M_n$.
  It is clear that for every $m$ and  $n$ in $\N$,
$\E_m\E_n=\E_n\E_m =\E_{\min(n,m)}$.

The following definition isolates the main topic of this paper.
\begin{definition}
A non-commutative martingale with respect to the filtration
$(\M_n)_{n=1}^\infty$ is a sequence $x=(x_n)_{n=1}^\infty$ in
$L^1(\M,\T)$ such that:
$$\E_n(x_{n+1})=x_n \qquad \text{for all}\ n\geq 1.$$
Similarly, if for all $n\geq 1$, $x_n$ is self-adjoint and
$\E_n(x_{n+1})\leq x_n$ (respectively, $\E_n(x_{n+1})\geq x_n $),
then the sequence $(x_n)_{n=1}^\infty$ is called a supermartingale
(respectively, submartingale).
\end{definition}
If additionally, $x \in L^p(\M,\T)$  for some $1<p<\infty$, then
$x$ is called a $L^p$-martingale. In this case, we set
\[
\Vert x \Vert_p :=\sup_{n\geq 1} \Vert x_n \Vert_p.
\]
If $\Vert x \Vert_p<\infty$, then $x$ is called a bounded
$L^p$-martingale. The difference sequence $dx=(dx_n)_{n=1}^\infty$
of  a martingale  $x=(x_n)_{n=1}^\infty$ is defined by
 $$dx_n = x_n - x_{n-1}$$
with the usual convention that $x_0=0$.


Recall that a subset $S$ of $L^1(\M,\T)$ is said to be uniformly
integrable if it is bounded and for every sequence of projections
$(p_n)_{n=1}^\infty $ with $p_n \downarrow_n 0$ (for the strong
operator topology), we have $\lim_{n \to \infty } \sup\{\Vert p_n
h p_n \Vert_1; h \in S \}=0$ (\cite{Ran10}). It can be easily
verified that a martingale $x=(x_n)_{n=1}^\infty$ in $L^1(\M,\T)$
is uniformly integrable if and only if there exists $x_\infty \in
L^1(\M,\T)$ such that $x_n=\E_n(x_\infty)$ for all $n \geq 1$. In
this case, the sequence $(x_n)_{n=1}^\infty$ converges to
$x_\infty$ in $L^1(\M,\T)$. In particular, if $1<p<\infty$, then
every bounded $L^p$-martingale is of the form
$(\E_n(x_\infty))_{n=1}^\infty$ for some $x_\infty \in
L^p(\M,\T)$.

For some concrete natural examples of non-commutative martingales,
we refer to \cite{PX} and the recent survey in this topic
\cite{X4}.

We will now describe square functions of non-commutative
martingales.  Following Pisier and Xu \cite{PX}, we will consider
the following row and column versions of
 square functions: for a martingale $x=(x_n)_{n=1}^\infty$, we denote by
$dx$ the difference sequence as defined above. For $N\geq 1$, set
\[
S_{C,N} (x) = \left(\sum^{N}_{k=1}|dx_{k}|^{2}\right)^{{1}/{2}} \
\ \mathrm{and}\ \ \ S_{R,N}(x) =
\left(\sum^{N}_{k=1}|dx_{k}^{*}|^{2}\right)^{{1}/{2}}.
\]
Let $E[0, \T({\bf 1}))$ be a rearrangement invariant  (quasi-)
Banach function space on the interval $[0, \T({\bf 1}))$. For any
finite sequence $a=(a_n)_{n \geq 1}$ in $E(\M,\T)$, set
\[
\|a\|_{E(\M;l^{2}_{C})} = \left\|\left(\sum_{n \geq
1}|a_{n}|^{2}\right)^{{1}/{2}}\right\|_{E(\M,\T)}, \ \
\|a\|_{E(\M; l^{2}_{R})} = \left\| \left(\sum_{n \geq 1}
|a_{n}^{*}|^{2}\right)^{{1}/{2}}\right\|_{E(\M,\T)}.
\]

The difference sequence $dx$ belongs to $E(\M;l^{2}_{C})$
(respectively,  $E(\M; l^{2}_{R})$)  if and only if  the sequence
$(S_{C,n}(x))_{n=1}^\infty$ (respectively,
$(S_{R,n}(x))_{n=1}^\infty$) is a bounded  sequence in $E(\M,\T)$.
In this case,  the limit $S_{C}(x) =
(\sum^{\infty}_{k=1}|dx_{k}|^{2})^{{1}/{2}}$ (respectively,
$S_{R}(x) = (\sum^{\infty}_{k=1}|dx_{k}^{*}|^{2})^{{1}/{2}}$) is
an element of $E(\M,\T)$. These two versions of square functions
are very crucial in the subsequent sections.

\section{Main Results: The Finite Case}

In this section, we assume that $\M$ is a finite von Neumann
algebra and $\T$ is normalized  normal faithful trace on $\M$.

We will retain all notations introduced in the previous two
sections. In particular, all adapted sequences are understood to
be with respect to a fixed filtration of von Neumann subalgebras
of $\M$. The principal result of this paper is Theorem~\ref{main}
below. It answers the problem raised in \cite{Ran16}.

\begin{theorem}\label{main}
There is an absolute constant $K>0$ such that if
$x=(x_n)^\infty_{n=1}$ is a $L^2$-bounded martingale, then there
exist  two sequences $y=(y_n)^\infty_{n=1}$ and
$z=(z_n)^\infty_{n=1}$ such that:
\begin{itemize}
\item[(i)] $(y_n)_{n=1}^\infty$ and $(z_n)_{n=1}^\infty$
are $L^2$-bounded martingales;
\item[(ii)] for every $n \geq 1$, $x_n=y_n +z_n$;
\item[(iii)] $\left\| dy\right\|_{L^{2}(\M; l^2_C)}
+ \left\| dz \right \|_{L^{2}(\M; l^2_R)} \leq 2 \| x \|_2$;
\item[(iv)] $\left\| dy\right\|_{L^{1, \infty}(\M; l^2_C)}
+ \left\| dz \right \|_{L^{1,\infty}(\M; l^2_R)} \leq K \| x
\|_1$.
\end{itemize}
\end{theorem}

As in the martingale transforms, our approach depends very heavily
on a non-commutative version of the classical Doob  weak type
$(1,1)$ maximal inequality, due to Cuculescu \cite{Cuc} (which we
will recall below) . As noted in \cite{Ran16}, the general case
can be deduced easily from the special case of positive
martingale. Hence,  without loss of generality, we can and do
assume that the martingale $x=(x_n)_{n=1}^\infty$ is a positive
martingale and $\Vert x\Vert_1=1$.

We will divide the proof into two parts. In the first part, we
will provide a detailed description of  the concrete decomposition
of the martingale $(x_n)_{n=1}^\infty$ and point out that $(i)$,
$(ii)$, and $(iii)$ are easily verified from the construction. In
the second part, we will show that the decomposition satisfies the
conclusion $(iv)$ of the theorem.

\medskip

\noindent \textsc{$\bullet$ Construction of the Martingales
${(y_n)_{n=1}^\infty}$
 and ${(z_n)_{n=1}^\infty}$.}

We start with  the proposition (due to Cuculescu \cite{Cuc}) below
which can be viewed as a substitute for the classical weak type
$(1,1)$ boundedness of maximal functions. We will state a version
that incorporates the different properties that we need in the
sequel.
 A short proof of the form stated below can be
 found in \cite{Ran15}.

\begin{proposition}[\cite{Cuc}]\label{maximal}
For every $\lambda >0$, there exists a sequence of decreasing
projections $(q_{n}^{(\lambda)})^{\infty}_{n=1}$ in $\M$ with:
\begin{itemize}
\item[(a)] for every $n \geq 1$, $q_{n}^{(\lambda)} \in \M_{n}$;
\item[(b)] $q_{n}^{(\lambda)}=\ch_{(0,\lambda]}(q_{n-1}^{(\lambda)} x_{n}
q_{n-1}^{(\lambda)})$. In particular,      $q_{n}^{(\lambda)}$
commutes with $q_{n-1}^{(\lambda)} x_{n} q_{n-1}^{(\lambda)}$;
\item[(c)]
 $q_{n}^{(\lambda)} x_{n} q_{n}^{(\lambda)} \leq \lambda
q_{n}^{(\lambda)}$;
\item[(d)]
if we  set $q^{(\lambda)} = \bigwedge_{n=1}^{\infty}
q_{n}^{(\lambda)}$ then $\T ({\bf 1}-q^{(\lambda)}) \leq
\lambda^{-1}$.
\end{itemize}
\end{proposition}

We consider collections of sequences of pairwise disjoint
projections as follows: for $n\geq 1$, set
\begin{equation}
\begin{cases}
p_{0,n} &:=\displaystyle{ \bigwedge^\infty_{k=0}q_{n}^{(2^k)}},\  \text{and} \\
p_{i,n} &:=\displaystyle{ \bigwedge^\infty_{k=i}q_{n}^{(2^k)} -
\bigwedge^\infty_{k=i-1} q_{n}^{(2^k)} \quad \text{for $i\geq
1$}}.
\end{cases}
\end{equation}
Similarly,
\begin{equation}
\begin{cases}
p_{0} &:= \displaystyle{\bigwedge^\infty_{k=0}q^{(2^k)}}, \ \text{and} \\
p_{i} &:= \displaystyle{\bigwedge^\infty_{k=i}q^{(2^k)} -
\bigwedge\limits^\infty_{k=i-1} q^{(2^k)} \quad \text{for $i\geq
1$}.}
\end{cases}
\end{equation}

Useful properties of  the  sequences $(p_{i,n})^\infty_{i=0}$ and
$(p_{i})^\infty_{i=0}$, that are relevant for our proof, are
collected in the following proposition whose verification is
straightforward and therefore is left to the reader.

\begin{proposition}\label{disjointmax} For $n\geq 1$, the sequence of projections
 $(p_{i,n})^\infty_{i=0}$ (respectively, $(p_{i})^\infty_{i=0}$)
 are pairwise disjoint
with the following properties:
\begin{itemize}
\item[(a)] For every $n\geq 1$ and $i\geq 0$, $p_{i,n} \in \M_n$;
\item[(b)] $\sum^\infty_{i=0} p_{i,n}={\bf 1}$ and $\sum^\infty_{i=0} p_{i}={\bf 1}$ (for the strong operator topology);
\item[(c)] for every $n_0 \geq 1$,  $\sum^{n_0}_{i=0}p_{i,n} \leq q_{n}^{(2^{n_0})}$ and
$\sum^{n_0}_{i=0}p_{i} \leq q^{(2^{n_0})}$.
\end{itemize}
\end{proposition}

 Since
$\sum^\infty_{i=0}p_{i,n}= {\bf 1}$, we have that
$a=\sum^\infty_{j=0} \sum^\infty_{i=0} p_{i,n} a p_{j,n}$ for all
$a \in  L^2 (\M, \T)$ so clearly, $a=\sum^\infty_{j=0}
\sum^\infty_{i\leq j} p_{i,n} a p_{j,n} + \sum^\infty_{j=0}
\sum^\infty_{i>j} p_{i,n} a p_{j,n}$. Our construction is based in
this simple fact.

 Define the sequences $y=(y_n)_{n=1}^\infty$ and
$z=(z_n)_{n=1}^\infty$ as follows:
\begin{equation}\label{mainequation}
\begin{cases}dy_1 &:=\displaystyle{ \sum^\infty_{j=0}\sum_{i \leq j} p_{i,1} dx_1 p_{j,1}};  \\
dy_n &:=\displaystyle{ \sum^\infty_{j=0}\sum_{i \leq j} p_{i,n-1} dx_n p_{j,n-1}}  \quad \text{for $n\geq 2$};\\
 dz_1 &:= \displaystyle{\sum^\infty_{j=0}\sum_{i > j}
p_{i,1} dx_1 p_{j,1}};\\
  dz_n &:= \displaystyle{\sum^\infty_{j=0}\sum_{i > j}
p_{i,n-1} dx_n p_{j,n-1}} \quad \text{for $n\geq 2$}.
\end{cases}
\end{equation}

Triangular truncations were also the main tools for the
construction of the decomposition used in \cite{Ran16}. The new
adjustment we need, in order to achieve the decomposition into two
martingales, is the use of  sequence of mutually disjoint
projections $(p_{i,n-1})_{i=0}^\infty$ from $\M_{n-1}$ (when
$n\geq 2$) instead of $(p_{i,n})_{i=0}^\infty$ used in
\cite{Ran16}. This was already clear since \cite{Ran16} but we
were unable to verify Proposition~\ref{lemma3} below  at that
time.

Clearly, $dx_n = dy_n + dz_n$ for every $n \geq 1$, therefore
$x_n=y_n +z_n$ for every $n\geq 1$, hence $(ii)$ is verified.

 Let $n\geq 2$.
Since $(p_{i,n-1})_{i=0}^\infty$ are  mutually disjoint
projections in $\M_{n-1}$,  and  triangular truncations are
orthogonal projections in $L^2(\M,\T)$, for every  $a \in
L^2(\M,\T)$,
\[\sum_{j=0}^\infty\sum_{i\leq j}
p_{i,n-1}ap_{j,n-1}= \lim_{k\to \infty}\sum_{j=0}^k\sum_{i\leq j}
p_{i,n-1}ap_{j,n-1}.\]
We deduce that for every $a \in
L^2(\M,\T)$,
\[
\E_{n-1}\left(\sum_{j=0}^\infty\sum_{i\leq j}
p_{i,n-1}ap_{j,n-1}\right)= \sum_{j=0}^\infty\sum_{i\leq j}
p_{i,n-1}\E_{n-1}(a)p_{j,n-1}.
\]
In particular, as $dx_n \in L^2(\M_n,\T_n)$ and
$\E_{n-1}(dx_n)=0$,  it follows that $dy_n \in L^2(\M_n,\T_n)$ and
$\E_{n-1}(dy_n)=0$. A similar remark can be made for
$\E_{n-1}(dz_n)$. Hence,
 $(dy_n)_{n=1}^\infty$ and $(dz_n)_{n=1}^\infty$ are
martingale difference sequences, which verifies $(i)$.

To verify   $(iii)$, it is  enough to note from the boundedness of
the triangular truncations in $L^2(\M,\T)$ that
\[
\sum_{n=1}^\infty \|dy_n\|_2^2 \leq \sum_{n=1}^\infty \|dx_n\|_2^2
=\|x\|_{2}^2.
\]
Noting that a similar inequality is also valid for
$\sum_{n=1}^\infty \|dy\|_2^2$, the inequality $(iii)$ follows.
Thus the items $(i)$, $(ii)$, and $(iii)$ of Theorem~\ref{main}
are verified.

\medskip

\noindent \textsc{$\bullet$ Proof of the Weak-Type
Inequality~(iv).}

In order to prove Theorem~\ref{main}(iv), we will make several
reductions. First, we remark   that
 from the construction of the martingales $(y_n)^\infty_{n=1}$ and
$(z_n)^\infty_{n=1}$, the square functions $(\sum^\infty_{n=1}
|dz_n^*|^2)^{1/2}$ and $(\sum^\infty_{n=1} |dy_n|^2)^{1/2}$ have
the same form. This can easily be seen  from the following lemma
whose verification is just a notational adjustment of
\cite[Lemma~2.2]{Ran16} and is left to the interested reader.

\begin{lemma}\label{lemma1}
For the sequences defined above, we have:
\begin{itemize}
\item[(a)] $|dy_1|^2 =\sum^\infty_{l=0} \sum^\infty_{j=0} \sum_{i\leq \min(l,j)} p_{l,1} dx_1 p_{i,1} dx_1 p_{j,1}$;
\item[(b)] $|dy_n|^2 =\sum^\infty_{l=0} \sum^\infty_{j=0} \sum_{i\leq \min(l,j)} p_{l,n-1} dx_n p_{i,n-1} dx_n
p_{j,n-1}$ for $n \geq 2$;
\item[(c)]$|dz_1^*|^2
=\sum^\infty_{l=1}\sum^\infty_{j=1}\sum_{i<\min(l,j)} p_{l,1} dx_1
p_{i,1} dx_1 p_{j,1}$;
\item[(d)] $|dz_n^*|^2
=\sum^\infty_{l=1}\sum^\infty_{j=1}\sum_{i<\min(l,j)} p_{l,n-1}
dx_n p_{i,n-1} dx_n p_{j,n-1}$ for $n\geq 2$,
\end{itemize}
where the sums are taken in the measure topology.
\end{lemma}

From the preceding  lemma, we only have to show that there is  an
absolute constant $C_1>0$ such that:
\begin{equation}
\left\| dy\right\|_{L^{1, \infty}(\M; l^2_C)} \leq C_1.
\end{equation}

 According to the definition of the quasi-norm $\|\cdot
\|_{1,\infty}$, this is equivalent   to show the existence of a
numerical constant $C_1>0$ such that for every $\lambda> 0$,
\begin{equation}\label{maininequality}
\T(\ch_{(\lambda, \infty)} (S_C(y))) \leq {C_1}\lambda^{-1}.
\end{equation}

The proof basically follows  the steps used in  \cite{Ran16,Ran15}
but some non-trivial adjustments had to be made.

\noindent  $\diamondsuit$  First,
 we consider the particular case:
{\it $\lambda = 2^{n_0}$ for some $n_0 \geq 0$.}

 The proof  of this case consists of three fundamental steps that will be highlighted in
 three separate propositions.

To avoid dealing with convergence, we will show that there is an
absolute constant $C_0>0$ such that for every $N \geq 1$,
\begin{equation}
\T(\ch_{(\lambda, \infty)} (S_{C,N}(y))) \leq {C_0} 2^{-n_0}.
\end{equation}

Throughout the proof, $N\geq 1$ is fixed.  We will reduce first to
the case of difference sequence of a bounded sequence in
$L^\infty(\M,\T)$. For notational purpose,
 we will simply write, throughout the proof,  $(q_n)_{n=1}^\infty$ (respectively, $q$)
for the projections $(q_n^{(2^{n_0})})_{n=1}^\infty$
(respectively, $q^{(2^{n_0})}$).  Consider the projection
\begin{equation}
w_{n_0}=\sum_{i=0}^{n_0} p_i= \bigwedge_{k=n_0}^\infty q^{(2^k)},
\end{equation}
and the operator
\begin{equation}
\gamma= |\sum_{j=0}^{n_0} \sum_{i\leq j} p_{i,1} dx_1 p_{j,1}|^2 +
\sum_{n=2}^{N}|\sum_{j=0}^{n_0} \sum_{i\leq j} p_{i,n-1} dx_n
p_{j,n-1}|^2. \end{equation}

The first step is to reduce the inequality from $S_{C,N}(y)$ to
$\gamma$. The significance of such reduction is the fact that the
triangular truncations used in $\gamma$ are formed from
 collections of finitely many projections.

\begin{varprop}\label{lemma2} For every $\alpha \in (0,1)$ and every $\beta \in
(0,1)$,
\begin{equation*}
\T\left(\ch_{(2^{n_0},\infty)}(S_{C,N}(y)) \right) \leq
\alpha^{-1} \T\left(\ch_{(\beta 4^{n_0},\infty)}(\gamma)
 \right) + 2(1-\alpha)^{-1} 2^{-n_{0}}.
\end{equation*}
\end{varprop}
\begin{proof}
Let $S = S_{C,N}(y)^2=\sum^N_{n=1} |dy_n|^2$.  Write $ S^{1/2} =
S^{1/2} w_{n_{0}} + S^{1/2}({\bf 1}-w_{n_{0}})$ and apply
Lemma~\ref{splitting} to get
\begin{equation*}
\begin{split}
\T\left(\ch_{(2^{n_0},\infty)}(S^{1/2})\right)&\leq \alpha^{-1}
\T\left(\ch_{(\sqrt{\beta}2^{n_0},\infty)}(|S^{1/2}w_{n_0}|)\right)\\
&\ + (1-\alpha)^{-1}
\T\left(\ch_{((1-\sqrt{\beta})2^{n_0},\infty)}(|S^{1/2}({\bf
1}-w_{n_0})|)\right).
\end{split}
\end{equation*}
 Since $\ch_{((1-\sqrt{\beta})2^{n_{0}}, \infty)} (|S^{1/2}({\bf
1}-w_{n_{0}})|)$ is a subprojection of ${\bf 1}-w_{n_{0}}$, it
follows that
\[
\T\left(\ch_{(2^{n_{0}},\infty)}(S^{ 1/2}) \right) \leq
\alpha^{-1}
\T\left(\ch_{(\sqrt{\beta}2^{n_0},\infty)}(|S^{1/2}w_{n_0}|)\right)
 + (1-\alpha)^{-1} \T({\bf 1}-w_{n_{0}}).
\]
Note that $w_{n_{0}} = \sum^{n_{0}}_{i=0} p_i =
\bigwedge^\infty_{k=n_{0}} q^{(2^{k})}$ so ${\bf 1}-w_{n_{0}}=
\bigvee^\infty_{k=n_{0}}({\bf 1}-q^{(2^k)})$. By
Proposition~\ref{maximal}(d), $ \T({\bf 1}-w_{n_{0}}) \leq
\sum^\infty_{k=n_{0}} \T({\bf 1}-q^{(2^k)}) \leq
\sum^\infty_{k=n_{0}} 2^{-k} = 2 .  2^{-n_{0}}$. Combining with
the previous  estimate, we conclude
\begin{equation*}
\begin{split}
\T\left(\ch_{(2^{n_{0}},\infty)}(S^{ 1/2}) \right) &\leq
\alpha^{-1} \T\left(\ch_{(\sqrt{\beta}2^{n_{0}},\infty)}(|S^{1/2}
w_{n_{0}}|)\right) + 2(1-\alpha)^{-1} 2^{-n_0}\\
 &=
\alpha^{-1} \T\left(\ch_{(\beta4^{n_{0}},\infty)}(w_{n_{0}}S
w_{n_{0}})\right) + 2(1-\alpha)^{-1} 2^{-n_0}.
\end{split}
\end{equation*}
To complete the proof, we will show that    $w_{n_{0}}S
w_{n_{0}}=w_{n_{0}}\gamma w_{n_{0}}$.

In fact, from the form of $|dy_n|^2$ stated in Lemma~\ref{lemma1},
we can write:
\begin{equation*}
\begin{split}
w_{n_{0}} S_{C,N}^2(y) w_{n_{0}}&=\sum^\infty_{l=0}
\sum^\infty_{j=0} \sum_{i\leq min(l,j)} w_{n_{0}} p_{l,1} dx_1
p_{i,1} dx_1 p_{j,1} w_{n_{0}} \\
&\ + \sum^N_{n=2} \sum^\infty_{l=0} \sum^\infty_{j=0} \sum_{i\leq
min(l,j)} w_{n_{0}} p_{l,n-1} dx_n p_{i,n-1} dx_n p_{j,n-1}
w_{n_{0}}.
\end{split}
\end{equation*}
We claim that all the sums taken in the expression of $w_{n_{0}}
S_{C,N}^2(y) w_{n_{0}}$ above are finite sums. For this, we remark
that if $l>n_0$ and $s\geq 1$, then $w_{n_{0}} p_{l,s}
=p_{l,s}w_{n_{0}}=0$. In fact, as $p_{l,s}= \bigwedge_{k=l}^\infty
q_s^{(2^k)} - \bigwedge_{k=l-1}^\infty q_s^{(2^k)}$ and $q^{(2^k)}
\leq q_s^{(2^k)}$ for all $k \geq 1$, it is clear that
$w_{n_0}=\bigwedge_{k=n_0}^\infty q^{(2^k)}$ is a subprojection of
$\bigwedge_{k=l-1}^\infty q_s^{(2^k)}$ when $l>n_0$ and therefore
$w_{n_0} \perp p_{l,s}$.  With this observation, we can write:
\begin{equation*}
\begin{split}
w_{n_{0}} S_{C,N}^2(y) w_{n_{0}}&=\sum^{n_0}_{l=0}
\sum^{n_0}_{j=0} \sum_{i\leq min(l,j)} w_{n_{0}} p_{l,1} dx_1
p_{i,1} dx_1 p_{j,1} w_{n_{0}} \\
&\ + \sum^N_{n=2} \sum^{n_0}_{l=0} \sum^{n_0}_{j=0} \sum_{i\leq
min(l,j)} w_{n_{0}} p_{l,n-1} dx_n p_{i,n-1} dx_n p_{j,n-1}
w_{n_{0}}\\
&=w_{n_{0}}\gamma w_{n_{0}}.
\end{split}
\end{equation*}
We conclude the proof by noting  that
\begin{equation*}
\begin{split}
\T\left(\ch_{(\beta 4^{n_0},\infty)}(w_{n_0}\gamma
w_{n_0})\right)&=\int_0^1 \ch_{(\beta
4^{n_0},\infty)}\left(\mu_t(w_{n_0}\gamma w_{n_0})\right)\ dt \\
&\leq \int_0^1 \ch_{(\beta
4^{n_0},\infty)}\left(\mu_t(\gamma)\right)\ dt \\
&= \T\left(\ch_{(\beta 4^{n_0},\infty)}(\gamma)\right).
\end{split}
\end{equation*}
The proof is complete.
\end{proof}

The next step is to   estimate
$\T(\ch_{(\beta4^{n_0},\infty)}(\gamma))$ using the $L^2$-norm of
square function of a supermartingale. This is the most significant
adjustment of the proof.

\begin{varprop}\label{lemma3}
The sequence $(q_n x_n q_n)^\infty_{n=1}$ is a supermartingale  in
$L^2(\M,\T)$  and if we set
$K_1=4\alpha^{-1}(1-\alpha)^{-1}(1-\beta)^{-1}\beta^{-2}
+10\sqrt{2}(1-\alpha)^{-2}(1-\beta)^{-1}(\sqrt{\beta})^{-1}$ for
$\alpha \in (0,1)$ and $\beta\in(0,1)$, then the following
inequality holds:
\begin{equation*}
\begin{split}
 \T\left( \ch_{(\beta 4^{n_0}, \infty)}(\gamma)\right) &\leq 2\alpha^{-1}\beta^{-2} 4^{-n_0}(\|q_1 x_1 q_1\|^2_2 +
\sum^N_{n=2} \|q_n x_n q_n-q_{n-1} x_{n-1} q_{n-1} \|^2_2) \\
&\ + K_12^{-n_0}.
\end{split}
\end{equation*}
\end{varprop}

\begin{proof}
The fact that  the sequence $(q_nx_nq_n)_{n=1}^\infty$ is a
supermartingale was already noted and proved in
\cite[Lemma~3.3]{Ran15} (see also, \cite[Lemma~2.4]{Ran16}) so
there is no need to repeat it here. To prove the estimate on
$\T\left( \ch_{(\beta 4^{n_0}, \infty)}(\gamma)\right) $, we note
that since for $0\leq j\leq n_0$ and $2\leq  n \leq N$, $p_{j,n-1}
\leq q_{n-1}$ and $q_n \leq q_{n-1}$, we have $p_{j,n-1}=
q_np_{j,n-1} + (q_{n-1}-q_n)p_{j,n-1}=p_{j,n-1}q_n +
p_{j,n-1}(q_{n-1}-q_n)$. We can decompose $\gamma$ as follows:
\begin{equation*}
\begin{split}
\gamma &= |\sum_{j=0}^{n_0} \sum_{i\leq j} p_{i,1} dx_1 p_{j,1}|^2
+ \sum_{n=2}^{N}|\sum_{j=0}^{n_0} \sum_{i\leq j} p_{i,n-1} dx_n
p_{j,n-1}|^2 \\
&=|\sum_{j=0}^{n_0} \sum_{i\leq j} p_{i,1} dx_1 p_{j,1}|^2 +
\sum_{n=2}^{N}|\sum_{j=0}^{n_0} \sum_{i\leq j} p_{i,n-1} dx_n(q_n)
p_{j,n-1} \\
&\ \ + \sum_{j=0}^{n_0} \sum_{i\leq j} p_{i,n-1} dx_n(q_{n-1}-q_n)
p_{j,n-1}|^2.
\end{split}
\end{equation*}
From the elementary inequality $|a +b|^2 \leq 2|a|^2 +2|b|^2$ for
any operators $a$ and $b$, we have,
\begin{equation*}
\begin{split}
\gamma &\leq |\sum_{j=0}^{n_0} \sum_{i\leq j} p_{i,1} dx_1
p_{j,1}|^2
 + 2\sum_{n=2}^{N}|\sum_{j=0}^{n_0} \sum_{i\leq j} p_{i,n-1}
dx_n(q_n) p_{j,n-1}|^2 \\
&\ \ + 2\sum_{n=2}^{N}|\sum_{j=0}^{n_0} \sum_{i\leq j} p_{i,n-1}
dx_n(q_{n-1}- q_n) p_{j,n-1}|^2.\\
\end{split}
\end{equation*}
The last of the three terms above  can be further decomposed to
get,
 \begin{equation*}
 \begin{split}
\gamma &\leq |\sum_{j=0}^{n_0} \sum_{i\leq j} p_{i,1} dx_1
p_{j,1}|^2
 + 2\sum_{n=2}^{N}|\sum_{j=0}^{n_0} \sum_{i\leq j} p_{i,n-1}
dx_n(q_n) p_{j,n-1}|^2 \\
&\ \ + 4\sum_{n=2}^{N}|\sum_{j=0}^{n_0} \sum_{i\leq j} p_{i,n-1}
q_ndx_n(q_{n-1}- q_n) p_{j,n-1}|^2 \\
&\ \ + 4\sum_{n=2}^{N}|\sum_{j=0}^{n_0} \sum_{i\leq j} p_{i,n-1}
(q_{n-1}-q_n)dx_n(q_{n-1}- q_n) p_{j,n-1}|^2.
\end{split}
\end{equation*}
Consider the following operators:
\begin{equation}\label{gamma123}
\begin{split}
\gamma_1 &:=|\sum_{j=0}^{n_0} \sum_{i\leq j} p_{i,1} dx_1
p_{j,1}|^2
 + 2\sum_{n=2}^{N}|\sum_{j=0}^{n_0} \sum_{i\leq j} p_{i,n-1}
dx_n(q_n) p_{j,n-1}|^2\\
\gamma_2 &:=4\sum_{n=2}^{N}|\sum_{j=0}^{n_0} \sum_{i\leq j}
p_{i,n-1} q_ndx_n(q_{n-1}- q_n) p_{j,n-1}|^2\\
\gamma_3 &:=4\sum_{n=2}^{N}|\sum_{j=0}^{n_0} \sum_{i\leq j}
p_{i,n-1} (q_{n-1}-q_n)dx_n(q_{n-1}- q_n) p_{j,n-1}|^2.
\end{split}
\end{equation}
Clearly, $\gamma\leq \gamma_1 +\gamma_2 +\gamma_3$.   The
splitting technique from Lemma~\ref{splitting} can be applied to
deduce that:
\begin{equation*}
\begin{split}
\T\left( \ch_{(\beta 4^{n_0}, \infty)}(\gamma)\right) &\leq
\alpha^{-1}\T\left( \ch_{(\beta^2 4^{n_0},
\infty)}(\gamma_1)\right) + (1-\alpha)^{-1}\T\left(
\ch_{((1-\beta)\beta
4^{n_0}, \infty)}(\gamma_2 +\gamma_3)\right)\\
&\leq \alpha^{-1}\T\left( \ch_{(\beta^2 4^{n_0},
\infty)}(\gamma_1)\right) + \alpha^{-1}(1-\alpha)^{-1}\T\left(
\ch_{((1-\beta)\beta^2 4^{n_0}, \infty)}(\gamma_2)\right) \\
&\ + (1-\alpha)^{-2}\T\left( \ch_{((1-\beta)^2\beta 4^{n_0},
\infty)}(\gamma_3)\right)\\
&= I +II + III.
\end{split}
\end{equation*}

We will estimate the quantities $I$, $II$, and $III$ in separate
three lemmas.

 \begin{lemma}\label{lemmaI}
 $I \leq 2\alpha^{-1}\beta^{-2}4^{-n_0}(\|q_1x_1q_1\|_2^2 +
 \sum^N_{n=2} \|q_n x_n q_n-q_{n-1} x_{n-1} q_{n-1} \|^2_2)$.
 \end{lemma}

 To prove this lemma, we remark first that
 \begin{equation*}
 \begin{split}
 I &=\alpha^{-1}\T\left( \ch_{(\beta^2 4^{n_0},
\infty)}(\gamma_1)\right)\\
&\leq  \alpha^{-1}\beta^{-2}4^{-n_0}(\|\sum^{n_0}_{j=0}
\sum_{i\leq j}  p_{i,1} dx_1 p_{j,1}\|_2^2 + 2\sum_{n=2}^N
\|\sum^{n_0}_{j=0} \sum_{i\leq j}  p_{i,n-1} dx_nq_n p_{j,n-1}
\|_2^2).
\end{split}
\end{equation*}
 Note that since  triangular
truncations are contractive in $L^2(\M,\T)$, the preceding
inequality yields:
\begin{equation*}
\begin{split}
I &\leq \alpha^{-1}\beta^{-2}4^{-n_0}(\|\sum^{n_0}_{j=0}
\sum_{i=0}^{n_0}  p_{i,1} dx_1 p_{j,1}\|_2^2 + 2 \sum_{n=2}^N
\|\sum^{n_0}_{j=0} \sum_{i=0}^{n_0} p_{i,n-1} dx_nq_n
p_{j,n-1}\|_2^2)\\
&\leq \alpha^{-1}\beta^{-2}4^{-n_0}(\|\sum^{n_0}_{j=0}
\sum_{i=0}^{n_0}  p_{i,1} dx_1 p_{j,1}\|_2^2 + 2 \sum_{n=2}^N \|
(\sum_{i=0}^{n_0} p_{i,n-1}) dx_nq_n \|_2^2).
\end{split}
\end{equation*}
Since for $j\geq 1$, $ \sum_{i=0}^{n_0}
p_{i,j}=\bigwedge_{k=n_0}^\infty q_j^{(2^k)}\leq q_j$
(Proposition~3.3(c)), we have
\begin{equation*}
\begin{split}
I &\leq \alpha^{-1}\beta^{-2}4^{-n_0}(\|q_1x_1q_1\|_2^2 +
2\sum_{n=2}^N \T\left(q_ndx_nq_{n-1}dx_n q_n\right))\\
&\leq 2\alpha^{-1}\beta^{-2}4^{-n_0}(\|q_1x_1q_1\|_2^2 +
\sum_{n=2}^N \T\left(q_ndx_nq_{n-1}dx_n q_n\right)).
\end{split}
\end{equation*}

To conclude the estimate on $I$, we will verify that for every
$n\geq 2$,
\[\T\left(q_n dx_n q_{n-1} dx_n q_n\right)\leq \|q_n x_n
q_n - q_{n-1}x_{n-1} q_{n-1} \|^2_2.\] This follows directly from
the facts (Proposition~\ref{maximal}) that $ q_n \leq q_{n-1}$ and
$q_n$ commutes with $q_{n-1}x_nq_{n-1}$. In fact,
\begin{equation*}
\begin{split}
 \T&\left(q_n
dx_n q_{n-1} dx_n q_{n}\right)
=\T\left(q_n (x_n-x_{n-1}) q_{n-1}(x_n-x_{n-1})q_n\right)\\
&=\T\left(q_n[q_{n-1}x_n q_{n-1}-
q_{n-1}x_{n-1}q_{n-1}][q_{n-1}x_n q_{n-1}-
q_{n-1}x_{n-1}q_{n-1}]q_n \right)\\
&\leq \T\left(q_n[q_{n}x_n q_{n}- q_{n-1}x_{n-1}q_{n-1}][q_{n}x_n
q_{n}- q_{n-1}x_{n-1}q_{n-1}]q_n \right)\\
&\leq \|q_nx_nq_n - q_{n-1}x_{n-1}q_{n-1}\|_2^2.
\end{split}
\end{equation*}
This shows that $I\leq
2\alpha^{-1}\beta^{-2}4^{-n_0}(\|q_1x_1q_1\|_2^2 + \sum_{n=2}^N
\|q_nx_nq_n - q_{n-1}x_{n-1}q_{n-1}\|_2^2)$ as stated in the
lemma.

\begin{lemma}\label{lemmaII}
 $II \leq  4\alpha^{-1}(1-\alpha)^{-1}(1-\beta)^{-1}\beta^{-2}2^{-n_0}$.
 \end{lemma}

To prove this estimate, recall that
$\gamma_2=4\sum_{n=2}^{N}|\sum_{j=0}^{n_0} \sum_{i\leq j}
p_{i,n-1}q_n dx_n(q_{n-1}- q_n) p_{j,n-1}|^2$. The proof rests
upon the following elementary but crucial  observation: for $2\leq
n\leq N$,
\begin{equation*}
q_ndx_n(q_{n-1}-q_n)=-q_nx_{n-1}(q_{n-1}-q_n).
\end{equation*}
Indeed, from the fact that $q_n$ commutes with $q_{n-1}x_n
q_{n-1}$, we have
\begin{equation*}
\begin{split}
q_nx_n(q_{n-1}-q_n)&=q_n(q_{n-1}x_nq_{n-1})(q_{n-1}-q_{n})\\
&=(q_{n-1}x_nq_{n-1})q_n(q_{n-1}-q_{n})=0.
\end{split}
\end{equation*}
We can now estimate $II$ as follows:
\begin{equation*}
\begin{split}
II &=\alpha^{-1}(1-\alpha)^{-1}\T\left( \ch_{((1-\beta)\beta^2
4^{n_0}, \infty)}(\gamma_2)\right)\\
&\leq 4\alpha^{-1}(1-\alpha)^{-1}(1-\beta)^{-1}\beta^{-2}4^{-n_0}
\sum_{n=2}^{N}\|\sum_{j=0}^{n_0} \sum_{i\leq j} p_{i,n-1}q_n
dx_n(q_{n-1}- q_n) p_{j,n-1}\|_2^2 \\
&\leq 4\alpha^{-1}(1-\alpha)^{-1}(1-\beta)^{-1}\beta^{-2}4^{-n_0}
\sum_{n=2}^{N}\|q_n
dx_n(q_{n-1}- q_n)\|_2^2 \\
 &= 4\alpha^{-1}(1-\alpha)^{-1}(1-\beta)^{-1}\beta^{-2}4^{-n_0}
\sum_{n=2}^N \|q_nx_{n-1}(q_{n-1}-q_n)\|_2^2\\
&=4\alpha^{-1}(1-\alpha)^{-1}(1-\beta)^{-1}\beta^{-2}4^{-n_0}
\sum_{n=2}^N \|q_n(q_{n-1}x_{n-1}q_{n-1})(q_{n-1}-q_n)\|_2^2\\
&\leq 4\alpha^{-1}(1-\alpha)^{-1}(1-\beta)^{-1}\beta^{-2}4^{-n_0}
\sum_{n=2}^N
\|q_n(q_{n-1}x_{n-1}q_{n-1})\|_\infty^2\|q_{n-1}-q_n\|_2^2.
\end{split}
\end{equation*}
By Proposition~\ref{maximal}(c), $\|q_{n-1}x_{n-1}q_{n-1}\|_\infty
\leq 2^{n_0}$ and therefore  we have
\begin{equation*}
\begin{split}
 II &\leq
4\alpha^{-1}(1-\alpha)^{-1}(1-\beta)^{-1}\beta^{-2}
\sum_{n=2}^N \T(q_{n-1}-q_n)\\
&=4\alpha^{-1}(1-\alpha)^{-1}(1-\beta)^{-1}\beta^{-2}
\T(q_{1}-q_N)\\
&\leq 4\alpha^{-1}(1-\alpha)^{-1}(1-\beta)^{-1}\beta^{-2}2^{-n_0}.
\end{split}
\end{equation*}
This completes the proof of the lemma.

\begin{lemma}\label{lemmaIII}
$III\leq
10\sqrt{2}(1-\alpha)^{-2}(1-\beta)^{-1}\beta^{-1}2^{-n_0}$.
\end{lemma}
The main tool for the proof is Proposition~\ref{truncation}. For
$2\leq n\leq N$, let $\cal{P}^{(n)}= \{p_{i,n-1}\}_{i=0}^{n_0}$.
Using the notation of Proposition~\ref{truncation}, $\gamma_3$ as
defined in (\ref{gamma123}) can be expressed as:
\begin{equation*}
\begin{split}
\gamma_3&=4\sum_{n=2}^{N}|\sum_{j=0}^{n_0} \sum_{i\leq j}
p_{i,n-1} (q_{n-1}-q_n)dx_n(q_{n-1}- q_n) p_{j,n-1}|^2\\
&= 4\sum_{n=2}^N
|\cal{T}^{\cal{P}^{(n)}}[(q_{n-1}-q_n)dx_n(q_{n-1}-q_n)]|^2\\
&= \sum_{n=2}^N
|\cal{T}^{\cal{P}^{(n)}}[2(q_{n-1}-q_n)dx_n(q_{n-1}-q_n)]|^2.
\end{split}
\end{equation*}
The crucial fact here is that for $2\leq n \leq N$,
\[
(q_{n-1}-q_n)dx_n(q_{n-1}-q_n)\geq 0.
\]
In fact, from the construction of the sequence of projection
$(q_n)_{n=1}^\infty$ from Proposition~\ref{maximal}(b), we have
$(q_{n-1}-q_n)x_n(q_{n-1}-q_n)\geq 2^{n_0}(q_{n-1}-q_n)$. On the
other hand, $q_{n-1}x_{n-1}q_{n-1}\leq 2^{n_0}q_{n-1}$ and
therefore, $(q_{n-1}-q_n)x_{n-1}(q_{n-1}-q_n)\leq
2^{n_0}(q_{n-1}-q_n)$. Hence,
$$(q_{n-1}-q_n)x_n(q_{n-1}-q_n)\geq
2^{n_0}(q_{n-1}-q_n)\geq (q_{n-1}-q_n)x_{n-1}(q_{n-1}-q_n)
$$
which shows that $(q_{n-1}-q_n)dx_n(q_{n-1}-q_n)\geq 0$. Therefore
Proposition~\ref{truncation} applies to $\gamma_3^{1/2}$. Hence,
we have the following estimates:
\begin{equation*}
\begin{split}
III &=(1-\alpha)^{-2}\T\left( \ch_{((1-\beta)^2\beta 4^{n_0},
\infty)}(\gamma_3)\right)\\
&=(1-\alpha)^{-2}\T\left( \ch_{((1-\beta)\sqrt{\beta} 2^{n_0},
\infty)}(\gamma_3^{1/2})\right)\\
&\leq
5\sqrt{2}(1-\alpha)^{-2}(1-\beta)^{-1}(\sqrt{\beta})^{-1}2^{-n_0}\sum_{n=2}^N
\|2(q_{n-1}-q_n)dx_n(q_{n-1}-q_n)\|_1\\
&\leq
10\sqrt{2}(1-\alpha)^{-2}(1-\beta)^{-1}(\sqrt{\beta})^{-1}2^{-n_0}\sum_{n=2}^N
\T\left((q_{n-1}-q_n)dx_n(q_{n-1}-q_n)\right).
\end{split}
\end{equation*}
Note that for every $2\leq n\leq N$ (using the fact that
$\E_{n-1}$ is $\T$-invariant),
\begin{equation*}
\begin{split}
\T\left((q_{n-1}-q_n)dx_n(q_{n-1}-q_n)\right)&=
\T\left((q_{n-1}-q_n)(x_n-x_{n-1})\right)\\
&=\T\left(q_{n-1}x_n -q_nx_n -q_{n-1}x_{n-1} +q_nx_{n-1}\right)\\
&=\T\left(q_{n-1}\E_{n-1}(x_n) -q_nx_n -q_{n-1}x_{n-1} +q_nx_{n-1}\right)\\
&=\T\left(-q_nx_nq_n + q_nx_{n-1}q_n\right)\\
&\leq \T\left(-q_nx_nq_n + q_{n-1}x_{n-1}q_{n-1}\right).
\end{split}
\end{equation*}
Taking the sum, we conclude that
\begin{equation*}
\begin{split}
III &\leq
10\sqrt{2}(1-\alpha)^{-2}(1-\beta)^{-1}(\sqrt{\beta})^{-1}2^{-n_0}\T\left(q_1x_1q_1
- q_{N}x_{N}q_{N}\right)\\
&=10\sqrt{2}(1-\alpha)^{-2}(1-\beta)^{-1}(\sqrt{\beta})^{-1}2^{-n_0}\T\left(
(q_1- q_{N})x_{N}\right)\\
&\leq
10\sqrt{2}(1-\alpha)^{-2}(1-\beta)^{-1}(\sqrt{\beta})^{-1}2^{-n_0}.
\end{split}
\end{equation*}
This completes the proof of the lemma.


The inequality in Proposition~B follows by  combining the
estimates on $I$, $II$, and $III$. The proof of the
 proposition is
complete.
\end{proof}

The last step is to estimate the $L^2$-norm of the square function
of the supermartingale $(q_nx_nq_n)_{n=1}^\infty$. This was
already achieved in \cite{Ran16} but we include below a much
shorter simplification that produces a better bound.

\begin{varprop}\label{lemma4} The square function of the
supermartingale from Proposition~\ref{lemma3} is $L^2$-bounded
with:
$$\left\| q_1x_1q_1 \right\|^{2}_{2}  +
\sum^N_{n=2} \|q_n x_n q_n-q_{n-1} x_{n-1} q_{n-1} \|^2_2\leq
  2^{n_0 +1}.$$
\end{varprop}
\begin{proof}
We will use the elementary identity $(a-b)^{2} = a^{2} - b^{2} +
b(b-a) + (b-a)b$ for self-adjoint operators. With $a=q_nx_nq_n$
and $b=q_{n-1}x_{n-1}q_{n-1}$, we have for every $n\geq 2$,
\begin{equation*}
\begin{split}
 \|q_nx_nq_n -q_{n-1}x_{n-1}q_{n-1}&\|_2^2
= \T\left((q_{n} x_{n} q_{n})^{2} - (q_{n-1} x_{n-1}
q_{n-1})^{2}\right) \\
&\quad + 2 \T\left(q_{n-1} x_{n-1} q_{n-1}[q_{n-1} x_{n-1} q_{n-1}
- q_{n} x_{n} q_{n}]\right) \\
 &=  \T\left((q_{n} x_{n} q_{n})^{2} -
(q_{n-1} x_{n-1} q_{n-1})^{2}\right)\\
 &\quad + 2 \T\left(q_{n-1}
x_{n-1} q_{n-1}[q_{n-1} x_{n-1} q_{n-1} - \E_{n-1}(q_{n} x_{n}
q_{n})]\right).
\end{split}
\end{equation*}
 By Proposition~\ref{maximal}(c),
$\left\|q_{n-1} x_{n-1} q_{n-1}\right\|_{\infty} \leq 2^{n_0}$.
Moreover, since the sequence $(q_nx_nq_n)_{n=1}^\infty$ is a
supermartingale,
 $q_{n-1} x_{n-1} q_{n-1} - \E_{n-1} (q_{n} x_{n} q_{n}) \geq 0$.
 Therefore, we get for every $n\geq 2$,
\begin{equation*}
\begin{split}
  \|q_nx_nq_n - x_{n-1}x_{n-1}q_{n-1}\|_2^2 &\leq
\T\left((q_{n} x_{n} q_{n})^{2} - (q_{n-1} x_{n-1}
q_{n-1})^{2}\right) \\
&\quad +  2^{n_0 +1} \T\left(q_{n-1} x_{n-1} q_{n-1} -
\E_{n-1}(q_{n} x_{n} q_{n})\right) \\
&=\T\left((q_{n} x_{n} q_{n})^{2} - (q_{n-1} x_{n-1}
q_{n-1})^{2}\right)  \\
&\quad  + 2^{n_0 +1} \T\left(q_{n-1} x_{n-1} q_{n-1} - q_{n} x_{n}
q_{n}\right).
\end{split}
\end{equation*}
Now,  we take the summation over $1\leq n \leq N$, we can conclude
that
\begin{equation*}
\begin{split}
\left\|q_1 x_1 q_1\right\|^{2}_2 &+ \sum_{n=2}^N \left\| q_nx_nq_n
- q_{n-1}x_{n-1}q_{n-1}\right\|^{2}_2 \\
 &\leq \left\|q_{1} x_{1} q_{1}\right\|^{2}_{2}
 + \sum^{N}_{n=2}
(\left\|q_{n} x_{n} q_{n}\right\|^{2}_{2} -
\left\|q_{n-1} x_{n-1} q_{n-1}\right\|^{2}_{2}) \\
&\ \ +  2^{n_0 +1}  \sum^{N}_{n=2} \T \left(q_{n-1} x_{n-1}
q_{n-1} - q_{n} x_{n}
q_{n}\right) \\
&= \left\|q_{1} x_{1} q_{1}\right\|^{2}_{2} +  (\left\|q_{N} x_{N}
q_{N}\right\|^{2}_{2} - \left\|q_{1} x_{1} q_{1}\right\|^{2}_{2})
+  2^{n_0 +1} \T\left(q_{1} x_{1} q_{1} -
q_{N} x_{N} q_{N}\right) \\
&= \left\|q_{N} x_{N} q_{N}\right\|^{2}_{2}  +
 2^{n_0 +1} \T \left((q_{1} - q_{N})x_{N}\right) \\
&\leq 2^{n_0 } \T\left(q_{N} x_{N}\right)  +
  2^{n_0 +1} \T \left((q_{1} - q_{N})x_{N} \right) \\
 &\leq   2^{n_0 +1}.
\end{split}
\end{equation*}
Thus the proof  is complete.
\end{proof}

We are now in a position to   conclude the proof of the weak-type
inequality (\ref{maininequality}) for the case $\lambda=
2^{n_{0}}$. This is accomplished by applying   successively
Proposition~\ref{lemma2}, Proposition~\ref{lemma3}, and
Proposition~\ref{lemma4} above. Indeed,
\begin{equation*}
\begin{split}
\T&\left(\ch_{(2^{n_0},\infty)}(S_{C,N}(y)) \right) \leq
\alpha^{-1}\T(\ch_{(\beta 4^{n_{0}}, \infty)}(\gamma))
+2(1-\alpha)^{-1} 2^{-n_{0}} \\
&\leq
\alpha^{-1}[2\alpha^{-1}\beta^{-2}4^{-n_{0}}(\|q_1x_1q_1\|_2^2 +
\sum_{n=1}^N \left\Vert q_nx_nq_n - q_{n-1}x_{n-1}q_{n-1}
\right\Vert_{2}^2) + K_1 2^{-n_0}]\\
 &\ \ +2(1-\alpha)^{-1} 2^{-n_{0}}\\
&\leq \alpha^{-1}\left[2\alpha^{-1}\beta^{-2}4^{-n_{0}}2^{n_0 +1} +
K_12^{-n_0}\right] +2(1-\alpha)^{-1} 2^{-n_{0}}\\
&=\left[4\alpha^{-2}\beta^{-2} + \alpha^{-1}K_1
+2(1-\alpha)^{-1}\right]2^{-n_{0}}.
\end{split}
\end{equation*}
If we set $C_0 := \inf\left\{4\alpha^{-2}\beta^{-2}
+\alpha^{-1}K_1 +2(1-\alpha)^{-1}; \alpha \in (0,1), \beta \in
(0,1) \right\}$ then
\begin{equation*}
\T\left(\ch_{(2^{n_0},\infty)}(S_{C,N}(y)) \right)\leq C_0
2^{-n_0}.
\end{equation*}
Hence taking the limit as $N$ tends to $\infty$,
inequality~(\ref{maininequality}) is verified for
$\lambda=2^{n_0}$.

$\diamondsuit$ Assume now the more general case that   $1 <
\lambda < \infty$.

\noindent Fix $n_0 \geq 0 $ such that $2^{n_{0}} < \lambda \leq
2^{n_{0}+1}$. We clearly have,
\begin{equation*}
\ch_{(\lambda,\infty)}(S_C(y)) \leq
\ch_{(2^{n_0},\infty)}(S_C(y)).
\end{equation*}
From the previous case, we can deduce,
\[
\T\left(\ch_{(\lambda,\infty)}(S_C(y))  \right) \leq C_0
2^{-n_{0}} = 2 C_0 2^{-(n_0+1)} \leq 2 C_0 \lambda^{-1}.
\]
Hence inequality (\ref{maininequality}) is verified for
$\lambda\geq 1$ with $C_1=2C_0$.


$\diamondsuit$ For the case  $0 < \lambda \leq 1$, we note that
since
 $\T$ is normalized, $\T\left( \ch_{(\lambda,\infty)}(S_C(y)) \right)\leq
 1$. In particular, $ \T\left( \ch_{(\lambda,\infty)}(S_C(y)) \right)
 \leq \lambda^{-1}$. Hence inequality~(\ref{maininequality}) is
 satisfied with a constant equals to $1$.


Combining the two cases $\lambda\geq 1$ and $\lambda<1$, we can
now conclude that
\[
\| dy \|_{L^{1, \infty}(\M;l_C^2)} \leq C_1.
\]
From the similarity of $|dy_n|^2$ and $|dz_n^*|^2$  demonstrated
in Lemma~\ref{lemma1}, we have
\[
\| dy \|_{L^{1, \infty}(\M;l_C^2)} + \| dz \|_{L^{1,
\infty}(\M;l_R^2)} \leq 2C_1=K_0.
\]
This completes the proof of Theorem~\ref{main}(iv) for the case of
normalized positive martingales.

The full generality  as stated in the theorem is obtained with
$K=8K_0$ by writing the martingale as linear combinations of four
positive martingales and normalization. Details are left to the
interested reader. \qed

\medskip

Recall that in general, triangular truncations are only of
weak-type $(1,1)$ (see for instance, \cite[Theorem~1.4]{DDPS2}).
The restriction to $L^2$-bounded martingales is needed in order to
verify that the sequences $(dy_n)_{n=1}^\infty$ and
$(dz_n)_{n=1}^\infty$ are martingale difference sequences. This
was possible from the boundedness of the triangular truncations.
If we replace $L^2(\M,\T)$ by any symmetric space  of measurable
operators  on which triangular truncations are bounded, then
$(dy_n)_{n=1}^\infty$ and $(dz_n)_{n=1}^\infty$, as constructed in
equation~(\ref{mainequation}),  are still martingale difference
sequences. However, the corresponding martingales may not be
bounded.
 Before proceeding, we need to recall the notion
of Boyd indices \cite[p.~130]{LT}. Let $E$ be a rearrangement
invariant Banach function space on $[0,1)$. For $s>0$, the
dilation operator $D_s: E \to E$ is defined by setting for any $f
\in E$,
\begin{equation*}
D_sf(t) = \begin{cases} f(t/s), \qquad &t\leq \min(1,s)\\
0, \quad &s<t<1\ (s<1).
\end{cases}
\end{equation*}
The {\it lower and upper Boyd indices } of $E$ are defined by
\begin{equation*}
\underline{\alpha}_E :=\lim_{s\to 0^{+}}\frac{\log\Vert
D_s\Vert}{\log s}, \qquad \overline{\alpha}_E :=\lim_{s\to
\infty}\frac{\log\Vert D_s\Vert}{\log s}.
\end{equation*}
It is well known that $0\leq \underline{\alpha}_E \leq
\overline{\alpha}_E \leq 1$ and if $E=L^p[0,1]$  for $1\leq p \leq
\infty$ then $\underline{\alpha}_E =\overline{\alpha}_E=1/p$. If
$0< \underline{\alpha}_E \leq \overline{\alpha}_E < 1$, we shall
say that $E$ has non-trivial Boyd indices. From
\cite[Theorem~3.3]{DDPS}, we can state:
\begin{theorem}\label{boyd}
There is an absolute positive constant $K>0$ such that if $E$ be a
rearrangement invariant Banach  function space on $[0,1)$ with the
Fatou property and   has non-trivial Boyd indices,
  and $x=(x_n)^\infty_{n=1}$ is a $L^1$-bounded martingale with
  $x_n \in E(\M,\T)$ for all $n\geq 1$,
 then there exist  two sequences $y=(y_n)^\infty_{n=1}$ and
$z=(z_n)^\infty_{n=1}$ such that:
\begin{itemize}
\item[($\alpha$)] for every $n \geq 1$, $x_n=y_n +z_n$;
\item[($\beta$)] $(y_n)_{n=1}^\infty$ and $(z_n)_{n=1}^\infty$
are martingales (not necessarily bounded);
\item[($\gamma$)] $\left\| dy\right\|_{L^{1, \infty}(\M; l^2_C)}
+ \left\| dz \right \|_{L^{1,\infty}(\M; l^2_R)} \leq K \| x
\|_1$.
\end{itemize}
\end{theorem}

Assume now that  $\M$ is hyperfinite and $(\M_n)_{n=1}^\infty$ is
a filtration consisting of finite dimensional von Neumann
subalgebras of $\M$, then the above restriction is no longer
needed.  In fact, in the case where the $\M_n$'s are finite
dimensional then for every $n\geq 2$, the mutually disjoint
sequence $(p_{i,n-1})_{i\geq 0}\subset \M_{n-1}$ (used in the
proof of Theorem~\ref{main}) is a finite sequence. Therefore, the
truncations used in the construction of the sequences
$(dy_n)_{n=1}^\infty$ and $(dz_n)_{n=1}^\infty$ are done with
finite sets of  mutually disjoint projections and consequently, is
bounded in $L^1(\M,\tau)$ (but not necessarily with uniform
bound). In this particular case, we can state the following result
as a complete
 non-commutative analogue of Theorem~0.1:
\begin{theorem}\label{main2}
There is an absolute constant $K$ such that if $\M$ is a finite
hyperfinite  von Neumann algebra and $(\M_n)_{n=1}^\infty$ is a
filtration in $\cal{M}$ consisting of finite dimensional von
Neumann subalgebras, then for every  $L^1$-bounded martingale
$x=(x_n)^\infty_{n=1}$, there exist two sequences
$y=(y_n)^\infty_{n=1}$ and $z=(z_n)^\infty_{n=1}$ such that:
\begin{itemize}
\item[($\alpha$)] for every $n \geq 1$, $x_n=y_n +z_n$;
\item[($\beta$)]  $ (y_n)_{n=1}^\infty$  and $(z_n)_{n=1}^\infty$
are $L^1$-martingales (not necessarily $L^1$-bounded);
\item[($\gamma$)] $\left\| dy\right\|_{L^{1, \infty}(\M; l^2_C)}
+ \left\| dz \right \|_{L^{1,\infty}(\M; l^2_R)} \leq K \| x
\|_1$.
\end{itemize}
\end{theorem}


\begin{remark}
 Theorem~\ref{main}  and Theorem~\ref{main2} can be extended to
square functions of non-commutative submartingales and
 non-commutative supermartingales. In this case, the
 decompositions are with submartingales (respectively,
 supermartingales). Details of such extension are done with just notational
 adjustments  of the  proof of
 \cite[Corollary~2.11]{Ran16} and are left to the interested reader.
\end{remark}

\section{Generalization to the semi-finite case}

In this section, we will consider the case where   $\M$ is  no
longer assumed to be finite. We can extend Theorem~\ref{main} to
the more general semi-finite case as follows:

\begin{theorem}\label{mainsemifinite}
There is an absolute constant $M>0$ such that if
$x=(x_n)^\infty_{n=1}$ is a martingale that is bounded in
$L^2(\M,\T) \cap L^1(\M,\T)$, then there exist two sequences
$v=(v_n)^\infty_{n=1}$ and $w=(w_n)^\infty_{n=1}$ such that:
\begin{itemize}
\item[(i)] $(v_n)_{n=1}^\infty$ and $(w_n)_{n=1}^\infty$
are $L^2$-bounded martingales;
\item[(ii)] for every $n \geq 1$, $x_n=v_n +w_n$;
\item[(iii)] $\left\| dv\right\|_{L^{2}(\M; l^2_C)}
+ \left\| dw \right \|_{L^{2}(\M; l^2_R)} \leq 2 \| x \|_2$;
\item[(iv)] $\left\| dv\right\|_{L^{1, \infty}(\M; l^2_C)}
+ \left\| dw \right \|_{L^{1,\infty}(\M; l^2_R)} \leq M \| x
\|_1$.
\end{itemize}
\end{theorem}

We will only outline the adjustments  needed for the proof   of
Theorem~\ref{main} to cover  the semi-finite case. For this, we
consider the case where $x=(x_n)_{n=1}^\infty$ is a positive
martingale and $\Vert x \Vert_1=1$ (the general case follows from
this case as noted at the end of the proof of Theorem~\ref{main}).
We will use the same notation as in the construction in Section~2
and Section~3. In particular, $(dy_n)_{n=1}^\infty$ and
$(dz_n)_{n=1}^\infty$ are martingale difference sequences  defined
as in (\ref{mainequation}).

We remark first that
 the fact that the trace $\T$  being  normalized, was used only to verify
 inequality~(\ref{maininequality}) when $0<\lambda <1$.
 That is,
 the proof of
 inequality~(\ref{maininequality}) when $\lambda\geq 1$ still
 applies for the semi-finite case.

As already noted in \cite{Ran16}, the only obstruction for proving
inequality~(\ref{maininequality}) without using the trace being
normalized  is the index $j=0$ in the definition of
$(dy_n)_{n=1}^\infty$. Indeed, if we set
\begin{equation}\label{firstsemifinite}
\begin{cases}ds_1 &:=\displaystyle{ \sum^\infty_{j=1}\sum_{i \leq j} p_{i,1} dx_1 p_{j,1}};  \\
ds_n &:=\displaystyle{ \sum^\infty_{j=1}\sum_{i \leq j} p_{i,n-1}
dx_n p_{j,n-1}}  \quad \text{for $n\geq 2$},
\end{cases}
\end{equation}
then   $dy_1 =p_{0,1}dx_1p_{0,1} + ds_1$ and  for $n\geq 2$, $dy_n
=p_{0,n-1}dx_np_{0,n-1} + ds_n$. Moreover, the sequence
$(s_n)_{n=1}^\infty$ is a $L^2$-bounded martingale and satisfies
the following weak-type inequality:

\begin{proposition}\label{lemmaS1} If $K$ is the positive constant
from Theorem~\ref{main}, then
 \[\left\Vert ds \right\Vert_{L^{1,\infty}(\M;l_C^2)}
+ \left\Vert dz \right\Vert_{L^{1,\infty}(\M;l_R^2)} \leq K.\]
\end{proposition}
\begin{proof}
As noted above, the fact that $\lambda \T\left( \ch_{(\lambda,
\infty)}(S_C(s))\right)\leq K$
  when $\lambda\geq 1$ is done exactly as in the
finite case.

For $0<\lambda<1$, we note that $|ds_1|$ is supported by the
projection ${\bf 1}-p_{0,1}$ and for $n\geq 2$,
 $|ds_n|$ is supported by the projection $({\bf
1}-p_{0,n-1})$. As $p_0 \leq p_{0,l}$ for every $l\geq 1$, it is
clear that $S_C(s)$ is supported by $({\bf 1}- p_0)$ and we claim
that $\T({\bf 1}- p_0) \leq 2$. This can be seen directly from the
definition of $p_0$. Indeed, $\T({\bf 1}-p_0)=\T({\bf
1}-\bigwedge_{k=0}^\infty q^{(2^k)})\leq \sum_{k=0}^\infty \T({\bf
1}- q^{(2^k)})\leq \sum_{k=0}^\infty 2^{-k}=2$. It follows that,
  for
$0<\lambda < 1$, $\lambda \T\left( \ch_{(\lambda,
\infty)}(S_C(\gamma))\right) \leq 2$. The same observation  on
support applies to $dz$ as well and thus the proof of the
proposition is complete.
\end{proof}

From Proposition~\ref{lemmaS1}, it is clear that we only need to
provide the \lq\lq right" decomposition of the martingale
difference sequence $(dy_n-ds_n)_{n=1}^\infty$. As in the
construction of $(dy_n)_{n=1}^\infty$ and  $(dz_n)_{n=1}^\infty$
in (\ref{mainequation}), we will decompose the projections $p_0$
and $p_{0,n}$'s into pairwise disjoint sequence of projections.
For $n\geq 1$ and $i\geq 0$, we set
\begin{equation}\label{defe}
\begin{cases}
e_{i,n} &:=\displaystyle{ \bigwedge^i_{k=0}(q_{n}^{(2^{-k})}
\wedge p_{0,n}) -
\bigwedge^{i+1}_{k=0} (q_{n}^{(2^{-k})}\wedge p_{0,n})};\\
e_{i} &:=\displaystyle{ \bigwedge^i_{k=0}(q^{(2^{-k})} \wedge
p_{0}) - \bigwedge^{i+1}_{k=0} (q^{(2^{-k})}\wedge p_{0})}.
\end{cases}
\end{equation}
Similarly,
\begin{equation}\label{defe2}
\begin{cases}
e_{-1,n} &:=\displaystyle{
\bigwedge^{\infty}_{k=0} (q_{n}^{(2^{-k})}\wedge p_{0,n})};\\
e_{-1} &:=\displaystyle{\bigwedge^{\infty}_{k=0}
(q^{(2^{-k})}\wedge p_{0})}.
\end{cases}
\end{equation}
\begin{remarks}\label{remarkS2}
 We have the following immediate properties:
\begin{itemize}
\item[(a)] For each $n \geq 1$, $(e_{i,n})_{i=-1}^\infty$ is a
sequence of disjoint projections and for $m\geq 1$, $\sum_{i=-1}^m
e_{i,n} = p_{0,n}- \bigwedge^{m+1}_{k=0}(q_{n}^{(2^{-k})} \wedge
p_{0,n}) + \bigwedge^{\infty}_{k=0}(q_{n}^{(2^{-k})} \wedge
p_{0,n}) $. In particular, $\sum_{i=-1}^\infty e_{i,n} = p_{0,n}$;
\item[(b)] For every $m \geq 1$,  $\sum_{i=m}^\infty e_{i,n} =
\bigwedge^{m+1}_{k=0}(q_{n}^{(2^{-k})} \wedge p_{0,n})-
\bigwedge^\infty_{k=0}(q_{n}^{(2^{-k})} \wedge p_{0,n})$. In
particular, $\sum_{k=m}^\infty e_{i,n} \leq q_n^{(2^{-m})}$.
\end{itemize}
\end{remarks}

It is clear from Remark~\ref{remarkS2} that for every $n \geq 2$,
$$
\sum_{j=-1}^\infty\sum_{i=-1}^\infty e_{i,n-1}dx_n e_{j,n-1} =
p_{0,n-1}dx_np_{0,n-1}.
$$
The decomposition of $(dy_n-ds_n)_{n=1}^\infty $ is done as in
(\ref{mainequation}), using the triangular truncation with respect
to the mutually disjoint sequence of  projections
$(e_{i,n-1})_{i=-1}^\infty$:
\begin{equation}\label{mainequationsemifinite}
\begin{cases}d\Xi_1 &:=\displaystyle{ \sum^\infty_{j=-1}\sum_{i \leq j} e_{i,1} dx_1 e_{j,1}};  \\
d\Xi_n &:=\displaystyle{ \sum^\infty_{j=-1}\sum_{i \leq j} e_{i,n-1} dx_n e_{j,n-1}}  \quad \text{for $n\geq 2$};\\
 d\Psi_1 &:= \displaystyle{\sum^\infty_{j=-1}\sum_{i > j}
e_{i,1} dx_1 e_{j,1}};\\
  d\Psi_n &:= \displaystyle{\sum^\infty_{j=-1}\sum_{i > j}
e_{i,n-1} dx_n e_{j,n-1}} \quad \text{for $n\geq 2$}.
\end{cases}
\end{equation}
Following the same line of argument used for the martingales
$(y_n)_{n=1}^\infty$ and $(z_n)_{n=1}^\infty$,  one can easily
verify that $(d\Xi_n)_{n=1}^\infty$ and $(d\Psi_n)_{n=1}^\infty$
are martingale difference sequences. Moreover,
$p_{0,1}dx_1p_{0,1}= d\Xi_1 + d\Psi_1$  and for every $n\geq 2$,
 \[p_{0,n-1}dx_np_{0,n-1}= d\Xi_n + d\Psi_n.\]
 Theorem~\ref{mainsemifinite} can be deduced from the following
 property of the  two  martingales $\Psi$ and $\Xi$.
\begin{proposition}\label{crucialsemifinite}
There is a numerical constant $C>0$ with:
\[
\left\| d\Psi\right\|_{L^{1, \infty}(\M; l^2_C)} + \left\|
d\Xi\right\|_{L^{1, \infty}(\M; l^2_R)}\leq C.
\]
\end{proposition}
Indeed,  if Proposition~\ref{crucialsemifinite} is verified, then
it is enough to set for $n\geq 1$, $v_n=\Psi_n +s_n$ and $w_n=
\Xi_n +z_n$ and Theorem~\ref{mainsemifinite}(iv) would follow
immediately from Proposition~\ref{lemmaS1} and
Proposition~\ref{crucialsemifinite}.

\smallskip

{\it Sketch of the proof of Proposition~\ref{crucialsemifinite}.}
First, we remark that as in Lemma~\ref{lemma1}, $|d\Psi_n|^2$ and
$|d\Xi_n^*|^2$ are of the same form. Therefore,
 as in the finite case,
  it is enough to verify that for every $0<\lambda
<\infty$,
\begin{equation}\label{maininequalitysemi}
\lambda\T\left(\ch_{(\lambda, \infty)}(S_C(\Psi))\right) \leq C.
\end{equation}

We will divide the proof into two cases.

$\diamondsuit$  \textsc{Case~1}: $\lambda \geq 1$. For  $N\geq 1$,
one can verify as in Lemma~3.4 that
\begin{equation*}
\begin{split}
 S_{C,N}^2(\Psi)  &= \sum^\infty_{l=-1}
\sum^\infty_{j=-1} \sum_{i\geq max(l,j)}  e_{l,1} dx_1 e_{i,1}
dx_1 e_{j,1}  \\
&\ +
 \sum^N_{n=2} \sum^\infty_{l=-1}
\sum^\infty_{j=-1} \sum_{i\geq max(l,j)}  e_{l,n-1} dx_n e_{i,n-1}
dx_n e_{j,n-1}.
\end{split}
\end{equation*}
As  $\sum_{j=-1}^\infty e_{j,n-1} = p_{0,n-1} \leq q_{n-1}^{(1)}$,
we can estimate $S_{C,N}^2(\Psi) \leq \bar{\gamma}_1 +
\bar{\gamma}_2 +\bar{\gamma}_3$ where the operators
$\bar{\gamma}_1$, $\bar{\gamma}_2$, and $\bar{\gamma}_3$ are
defined as follows:
\begin{equation*}
\begin{split}
\bar{\gamma}_1 &:=|\sum_{j=-1}^{\infty} \sum_{i> j} e_{i,1} dx_1
e_{j,1}|^2
 + 2\sum_{n=2}^{N}|\sum_{j=-1}^{\infty} \sum_{i> j} e_{i,n-1}
dx_n(q_n^{(1)}) e_{j,n-1}|^2\\
\bar{\gamma}_2 &:=4\sum_{n=2}^{N}|\sum_{j=-1}^{\infty} \sum_{i> j}
e_{i,n-1} q_n^{(1)}dx_n(q_{n-1}^{(1)}- q_n^{(1)}) e_{j,n-1}|^2\\
\bar{\gamma}_3 &:=4\sum_{n=2}^{N}|\sum_{j=-1}^{\infty} \sum_{i>j}
e_{i,n-1} (q_{n-1}^{(1)}-q_n^{(1)})dx_n(q_{n-1}^{(1)}- q_n^{(1)})
e_{j,n-1}|^2.
\end{split}
\end{equation*}
Using the splitting  technique from Lemma~\ref{splitting}, we can
deduce that for every $\alpha\in(0,1)$ and $\beta\in (0,1)$:
\begin{equation*}
\begin{split}
\T\left(\ch_{(\lambda, \infty)}(S_{C, N}(\Psi))\right)
&=\T\left(\ch_{(\lambda^2, \infty)}(S_{C,N}^2(\Psi))\right)\\
&\leq \alpha^{-1}\T\left(\ch_{(\beta\lambda^2,
\infty)}(\bar{\gamma}_1)\right) +
\alpha^{-1}(1-\alpha)^{-1}\T\left(\ch_{(\beta(1-\beta)\lambda^2,
\infty)}(\bar{\gamma}_2)\right) \\
&\ +
(1-\alpha)^{-2}\T\left(\ch_{(\beta^2\lambda^2, \infty)}(\bar{\gamma}_3)\right)\\
&= IV + V + VI.
\end{split}
\end{equation*}
We can estimate $IV$, $V$, and $VI$ separately  following the
proofs of  Lemma~\ref{lemmaI}, Lemma~\ref{lemmaII}, and
Lemma~\ref{lemmaIII} to  deduce that there are constants $A_1$ and
$A_2$ (depending only on $\alpha$ and $\beta$) such that
\begin{equation*}
\begin{split}
\T\left(\ch_{(\lambda, \infty)}(S_{C, N}(\Psi))\right) &\leq
A_1\lambda^{-2}(\|q_1^{(1)}x_1q_1^{(1)}\|_2^2 + \sum_{n=2}^\infty
\Vert q_{n}^{(1)}x_nq_{n}^{(1)} -
q_{n-1}^{(1)}x_{n-1}q_{n-1}^{(1)}\Vert_{2}^2) \\
&\ + A_2\lambda^{-1}.
\end{split}
\end{equation*}
Now, we can apply   Proposition~\ref{lemma4}  with $n_0=0$ to get
\begin{equation*}
 \|q_1^{(1)}x_1q_1^{(1)}\|_2^2 + \sum_{n=2}^\infty
\Vert q_{n}^{(1)}x_nq_{n}^{(1)} -
q_{n-1}^{(1)}x_{n-1}q_{n-1}^{(1)}\Vert_{2}^2\leq 2.
\end{equation*}
Combining  the last two inequalities, we conclude that
\begin{equation*}
\T\left(\ch_{(\lambda,\infty)}(S_C(\Psi)) \right) \leq
2A_1\lambda^{-2} + A_2\lambda^{-1}
\end{equation*}
and since $\lambda \geq 1$,  (\ref{maininequalitysemi}) follows.


$\diamondsuit$ \textsc{Case~2:} For $\lambda <1$, we will consider
the special case $\lambda=2^{-n_0}$ for some $n_0\geq 1$. Consider
the following projection
\begin{equation*}
f_{n_0} = \sum_{i=n_0 }^\infty e_{i},
\end{equation*}
and the operator
\begin{equation*}
\varphi_0= |\sum_{j=n_0}^\infty \sum_{i\geq j} e_{i,1} dx_1
e_{j,1}|^2 + \sum_{n=2}^{N}|\sum_{n_0}^\infty \sum_{i\geq j}
e_{i,n-1} dx_n e_{j,n-1}|^2.
\end{equation*}
It is easy to verify that $\T({\bf 1} -f_{n_0}) \leq 2^{n_0 +1}$.
Moreover, $f_{n_{0}} S_{C,N}^2(\Psi) f_{n_{0}}\leq f_{n_0}
\varphi_0 f_{n_0}$. This can be seen as follows: first, note that
\begin{equation*}
\begin{split}
f_{n_{0}} S_{C,N}^2(\Psi) f_{n_{0}}&=\sum^\infty_{l=0}
\sum^\infty_{j=0} \sum_{i\geq max(l,j)} f_{n_{0}} e_{l,1} dx_1
e_{i,1} dx_1 e_{j,1} f_{n_{0}}\\
&\ + \sum^N_{n=2} \sum^\infty_{l=0} \sum^\infty_{j=0} \sum_{i\geq
max(l,j)} f_{n_{0}} e_{l,n-1} dx_n e_{i,n-1} dx_n e_{j,n-1}
f_{n_{0}}.
\end{split}
\end{equation*}
Remark that if  $s\geq 1$ and  $l<n_0$, then $f_{n_{0}} e_{l,s}
=e_{l,s}f_{n_{0}}=0$. In fact, we note that  as $e_{l,s}=
\bigwedge_{k=0}^l (q_s^{(2^{-k})}\wedge p_{0,s}) -
\bigwedge_{k=0}^{l+1} (q_n^{(2^{-k})}\wedge p_{0,s})$,
$q^{(2^{-k})} \leq q_s^{(2^{-k})}$ for all $k \geq 1$ and $p_0
\leq p_{0,n}$, it is clear that $f_{n_0}=\bigwedge_{k=0}^{n_0 +1}
(q^{(2^{-k})}\wedge p_0) - \bigwedge_{k=0}^{\infty}
(q^{(2^{-k})}\wedge p_0)$ is a subprojection of
$\bigwedge_{k=0}^{l+1} (q_s^{(2^{-k})}\wedge p_{0,s})$ when
$l<n_0$ and by the definition of $e_{l,s}$, it follows that
$f_{n_0} \perp e_{l,s}$. Therefore,
\begin{equation*}
\begin{split}
f_{n_{0}} S_{C,N}^2(\Psi) f_{n_{0}}&=\sum^\infty_{l=n_0}
\sum^\infty_{j=n_0} \sum_{i\geq max(l,j)} f_{n_{0}} e_{l,1} dx_1
e_{i,1} dx_1 e_{j,1} f_{n_{0}}\\
&\ + \sum^N_{n=2} \sum^\infty_{l=n_0} \sum^\infty_{j=n_0}
\sum_{i\geq max(l,j)} f_{n_{0}} e_{l,n-1} dx_n e_{i,n-1} dx_n
e_{j,n-1} f_{n_{0}} \\
&=f_{n_{0}}\varphi_0f_{n_{0}}.
\end{split}
\end{equation*}
Write $S_{C,N}(\Psi)^2=f_{n_{0}} S_{C,N}^2(\Psi) f_{n_{0}} +
f_{n_{0}} S_{C,N}^2(\Psi)({\bf 1}- f_{n_{0}}) + ({\bf
1}-f_{n_{0}}) S_{C,N}^2(\Psi)$. Using the splitting techniques as
in Proposition~\ref{lemma2}, we can make the following reduction:
\begin{lemma}\label{lemmaS3}
For  $\alpha \in (0,1)$ and
 $\beta \in (0,1)$,
\begin{equation*}
\T\left(\ch_{(2^{-n_0},\infty)}(S_C(\Psi)) \right) \leq
\alpha^{-1}\T(\ch_{(\beta 4^{-n_0}, \infty)}(\varphi_0)) +
(1-\alpha)^{-1} 2^{n_{0}}.
\end{equation*}
\end{lemma}
From Remark~\ref{remarkS2}(iii), $ \sum_{i=n_0}^{\infty} e_{i,n-1}
\leq q_{n-1}^{(2^{-n_0})}$, and therefore as above, we have
$\varphi_0\leq \varphi_1 +\varphi_2 +\varphi_3$ where
\begin{equation*}
\begin{split}
 \varphi_1 &:=|\sum_{j=n_0}^{\infty} \sum_{i> j} e_{i,1} dx_1
e_{j,1}|^2
 + 2\sum_{n=2}^{N}|\sum_{j=n_0}^{\infty} \sum_{i> j} e_{i,n-1}
dx_n(q_n^{(2^{-n_0})}) e_{j,n-1}|^2\\
\varphi_2 &:=4\sum_{n=2}^{N}|\sum_{j=n_0}^{\infty} \sum_{i> j}
e_{i,n-1} q_n^{(2^{-n_0})}dx_n(q_{n-1}^{(2^{-n_0})}-
q_n^{(2^{-n_0})}) e_{j,n-1}|^2\\
 \varphi_3 &:=4\sum_{n=2}^{N}|\sum_{j=n_0}^{\infty}
\sum_{i>j} e_{i,n-1}
(q_{n-1}^{(2^{-n_0})}-q_n^{(2^{-n_0})})dx_n(q_{n-1}^{(2^{-n_0})}-
q_n^{(2^{-n_0})}) e_{j,n-1}|^2.
\end{split}
\end{equation*}

Using the splitting  technique from Lemma~\ref{splitting} a second
time, we can deduce that for $\alpha\in(0,1)$ and $\beta\in
(0,1)$:
\begin{equation*}
\begin{split}
\T\left(\ch_{(\beta 4^{-n_0}, \infty)}(\varphi_0)\right) &\leq
\alpha^{-1}\T\left(\ch_{(\beta^2 4^{-n_0},
\infty)}(\varphi_1)\right)\\
&\ +
\alpha^{-1}(1-\alpha)^{-1}\T\left(\ch_{(\beta^2(1-\beta)4^{-n_0},
\infty)}(\varphi_2)\right)\\
&\  +
(1-\alpha)^{-2}\T\left(\ch_{(\beta(1-\beta)^2 4^{-n_0}, \infty)}(\varphi_3)\right)\\
&= A + B + C.
\end{split}
\end{equation*}
We can estimate $A$, $B$, and $C$ as in Lemma~\ref{lemmaI},
Lemma~\ref{lemmaII}, and Lemma~\ref{lemmaIII}  to  deduce that
there are constants $B_1$ and $B_2$ (depending only on $\alpha$
and $\beta$) such that
\begin{equation*}
\begin{split}
\T\left(\ch_{(\lambda, \infty)}(S_{C, N}(\Psi))\right) &\leq B_1
4^{n_0} (\|q_1^{(2^{-n_0})}x_1q_1^{(2^{-n_0})}\|_2^2 \\
&\ + \sum_{n=2}^\infty \Vert
q_{n}^{(2^{-n_0})}x_nq_{n}^{(2^{-n_0})} -
q_{n-1}^{(2^{-n_0})}x_{n-1}q_{n-1}^{(2^{-n_0})}\Vert_{2}^2) + B_2
2^{n_0}.
\end{split}
\end{equation*}

Once the truncation in the preceding inequality is established, we
can proceed  exactly as in the proof of the finite case with
$2^{-n_0}$ instead of $2^{n_0}$. \qed

\medskip

We conclude this section with the obvious remark that
Theorem~\ref{boyd} and the extensions to supremartingales and
submartingales  still hold for the semi-finite case.

\section{Applications: Optimal orders of growth of the constants in
non-commutative Burkholder-Gundy inequalities}

Throughout  this section, we assume that $\M$ is finite and the
trace $\T$ is normalized.  We recall the definitions of martingale
Hardy spaces.  For $1 \leq p < \infty$, $\H^{p}_{C}(\M)$
 (respectively, $\H^{p}_{R}(\M)$) is defined as
the set of all $L^{p}$-martingales $x$ with respect to a
filtration $(\M_{n})_{n \geq 1}$ such that $dx \in L^{p}(\M;
l^{2}_{C})$ (respectively, $L^{p}(\M;l^{2}_{R})$), and set
\[
\|x\|_{\H^{p}_{C}(\M)} = \|dx\|_{L^{p}(\M;l^{2}_{C})} \ \
\text{and} \ \ \|x\|_{\H^{p}_{R}(\M)} =
\|dx\|_{L^{p}(\M;l^{2}_{R})}.
\]
Equipped with the previous norms, $\H^{p}_{C}(\M)$ and
$\H^{p}_{R}(\M)$ are Banach spaces. The Hardy space of
non-commutative martingale is defined as follows:  if $1 \leq p <
2$, \[\H^{p}(\M) = \H^{p}_{C}(\M) + \H^{p}_{R}(\M)\]
 equipped with
the norm
\[
\|x\|_{\H^{p}(\M)} = \inf\left\{\|y\|_{\H^{p}_{C}(\M)} +
\|z\|_{\H^{p}_{R}(\M)}\right\}\] where the infimum is taken over
all $y$ and $z$ with  $ x = y + z$, $y \in \H^{p}_{C}(\M)$, and $z
\in \H^{p}_{R}(\M)$;
 and if $2 \leq p < \infty $,
\[
\H^{p}(\M) = \H^{p}_{C}(\M) \cap \H^{p}_{R}(\M)\]
 equipped with the norm
\[\|x\|_{\H^{p}(\M)} = \max\left\{ \|x\|_{\H^{p}_{C}(\M)} , \
\|x\|_{\H^{p}_{R}(\M)}\right\}.\]

The aim of this section is to point out that using interpolation
techniques, Theorem~\ref{main} provides  a new proof of the
non-commutative Burkholder-Gundy inequalities.
\begin{theorem}[\cite{PX}]\label{BG}
Let $1 < p < \infty$.  Let $x = (x_{n})^{\infty}_{n=1}$ be an
$L^{p}$-martingale.  Then $x$ is bounded in $L^{p}(\M,\T)$ if and
only if $x$ belongs to $\H^{p}(\M)$.  If this is the case then,
\begin{equation*}
 \alpha^{-1}_{p}\|x\|_{\H^{p}(\M)} \leq \|x\|_{p} \leq
\beta_{p}\|x\|_{\H^{p}(\M)}. \tag {$BG_p$}
\end{equation*}
\end{theorem}
For Clifford martingales, some particular cases of
Theorem~\ref{BG} was also obtained in \cite{CK}.  Note that up
until now, both the original proof in \cite{PX} and the
alternative proof from \cite{Ran15} made use of  the
non-commutative Stein inequality (also proved in \cite{PX}) in
order to achieve the decomposition into two martingales, as
described in the definition of $ \|\cdot\|_{\H^{p}(\M)}$ for
$1<p<2$. Theorem~\ref{main} allows us to avoid the use of the
non-commutative Stein inequality. This approach, which is probably
more complex than the existing proofs, produces better estimates
of the constants involved. Indeed, it allows us to deduce the
optimal order of growth for the constant $\alpha_p$ (which is the
same as in the case of commutative case) when $p\to 1$. This
solves a question left open in \cite{JX2} (see the remark after
the main theorem of \cite{JX2}).

We will write $a_p\approx b_p$ as $p \to p_0$ to abbreviate the
statement that there are two absolute positive constants $K_1$ and
$K_2$ such that
$$
K_1\leq \frac{a_p}{b_p}\leq K_2 \ \text{for $p$ close to $p_0$}.
$$

The following theorem is the principal result of this section.
\begin{theorem}\label{alpha}
We have the following estimates for the best constants in
$(BG_p)$:
\begin{itemize}
\item[(i)] $\alpha_p\approx (p-1)^{-1}$ as $p \to 1$;
\item[(ii)] $\alpha_p\approx p$ as $p \to \infty$;
\item[(iii)] $\beta_p\approx 1$ as $p \to 1$;
\item[(iv)] $\beta_p\approx p$ as $p\to \infty$.
\end{itemize}
These are the optimal orders of growth of the constants $\alpha_p$
and $\beta_p$.
\end{theorem}

\begin{remarks}  (a) Compared with the commutative setting, the optimal
orders of $\beta_p$ are the same as its commutative counterpart.
However, $\alpha_p$ behaves differently. In the commutative case,
$\alpha_p\approx \sqrt{p}$ when $p\to \infty$, and
$\alpha_p\approx (p-1)^{-1}$ when $p \to 1$.

 (b) The only new result here is $(i)$.  The optimal orders as stated
in $(ii)$, $(iii)$, and $(iv)$  were obtained by combining results
from  \cite{JX}, \cite{JX2}, and \cite{Ran15}. We also note that
for the special case of even integers, $(ii)$ was established in
\cite{PS5} for more general sequences called $p$-orthogonal sums.
\end{remarks}

 We  will use the real interpolation, namely the $J$-method. We will
review the general theory of real interpolation. Our main
reference for interpolation is the book of Bergh and
L\"{o}fstr\"{o}m \cite{BL} and the recent survey \cite{KaltonSMS}.

A pair of (quasi)-Banach spaces $(E_0, E_1)$ is called a
compatible couple if they embed continuously into some topological
vector space $X$. This allows us to consider the spaces $E_0\cap
E_1$ and $E_0 + E_1$ equipped with $\|x\|_{E_0\cap
E_1}=\max\{\|x\|_{E_0}, \|x\|_{E_1}\}$, $\|x\|_{E_0 +
E_1}=\inf\{\|x_0\|_{E_0} + \|x_1\|_{E_1} : x=x_0 + x_1,\  x_0 \in
E_0,\ x_1 \in E_1\}$ respectively.

For a compatible couple $(E_0, E_1)$, we define for any $x \in E_0
\cap E_1$, and $t>0$,
\[
J(x,t;E_0, E_1)=\max\{\|x\|_{E_0}, t\|x\|_{E_1}\}.
\]
If the compatible couple is clear from the context, we will simply
write $J(x,t)$.

To avoid dealing with measurability, we will be working with the
{\it discrete version of the $J$-method} which we will now
describe: for $0<\theta<1$ and $1\leq p<\infty$, we denote by
$\lambda^{\theta,p}$ the space of all sequences
$(\alpha_\nu)_{\nu=-\infty}^{\infty}$  for which,
\[
\|(\alpha_\nu)\|_{\lambda^{\theta,p}}=\left\{\sum_{\nu\in\mathbb{Z}}
(2^{-\nu\theta}|\alpha_\nu|)^p \right\}^{\frac{1}{p}} <\infty.
\]
\begin{definition}
Let $(E_0,E_1)$ be a compatible couple and suppose that
$0<\theta<1$, and $1\leq p<\infty$. The interpolation space
$(E_0,E_1)_{\theta,p,\underline{J}}$ consists of elements $x \in
E_0 + E_1$ which admits a representation:
\begin{equation}\label{interpolation}
x=\sum_{\nu \in \mathbb{Z}} u_\nu, \quad \text{(convergence  in
$E_0 +E_1$)},
\end{equation}
with $u_\nu \in E_0 \cap E_1$ and such that
\[
\|x\|_{\theta,p,\underline{J}}=\inf\left\{\|\{J(u_\nu,
2^\nu)\}\|_{\lambda^{\theta,p}}\right\}<\infty,
\]
where the infimum  is taken over all representations of $x$ as in
(\ref{interpolation}).
\end{definition}

For general information on interpolations of non-commutative
spaces, we refer to \cite{DDP2} and \cite[p.~1466]{PX3}.

 \begin{proof}[Proof of Theorem~\ref{alpha}(i)] It is enough to
consider positive $L^2$-bounded martingale $x=(x_n)_{n=1}^\infty$.
Let $x_\infty\in L^2(\M,\T)$ such that $x_n=\E_n(x_\infty)$ for
every $n\geq 1$.  Let $1<p<2$ and $0<\theta<1$ with
$1/p=(1-\theta) +\theta/2$. For $\epsilon
>0$, fix $(u_\nu)_{\nu=-\infty}^\infty$ in $L^2(\M,\T)$ such that
\begin{equation*}
\begin{split}
x_\infty &=\sum_{\nu \in \mathbb{Z}} u_\nu \quad \text{and}\\
 \|x_\infty\|_{\theta,p;\underline{J}} +\epsilon &\geq
\|\{J(u_\nu, 2^\nu)\}\|_{\lambda^{\theta,p}},
\end{split}
\end{equation*}
where the $J$-functional is relative to the  interpolation couple
$(L^1(\M,\T), L^2(\M,\T))$.

For each $\nu \in \mathbb{Z}$, Theorem~\ref{main} guaranties the
existence of  an absolute constant $K>0$, and  two $L^2$-bounded
martingales $y^{(\nu)}$ and $z^{(\nu)}$ such that:
\begin{itemize}
\item[(a)] $\E_n(u_\nu)= y_n^{(\nu)} + z_n^{(\nu)}$ for all $n\geq 1$;
\item[(b)] $J(S_C(dy^{(\nu)}), t) \leq KJ(u_\nu,t)$ for every
$t>0$;
\item[(c)] $J(S_R(dz^{(\nu)}), t) \leq KJ(u_\nu,t)$ for every
$t>0$,
\end{itemize}
where the $J$-functionals in the left hand side of the
inequalities in (b) and (c)  above are taken relative to the
interpolation couple $(L^{1,\infty}(\M,\T), L^2(\M,\T))$. From
this, we can deduce that,
$$
\|\{J(S_C(dy^{(\nu)}), 2^\nu)\}\|_{\lambda^{\theta,p}} \leq
K(\|x_\infty\|_{\theta,p;\underline{J}} +\epsilon).
$$
Note   from the definition of the $J$-functionals that,
\begin{equation*}
\begin{split}
J&(S_C(dy^{(\nu)}),
2^\nu)=\max\left\{\|S_C(dy^{(\nu)})\|_{1,\infty}, 2^\nu\|S_C(dy^{(\nu)})\|_2\right\}\\
&=\max\left\{\|\sum_n dy_n^{(\nu)} \otimes e_{n,1}
\|_{L^{1,\infty}(\M\overline{\otimes} B(l^2)) }, 2^\nu\|\sum_n
dy_n^{(\nu)} \otimes e_{n,1}\|_{L^2(\M \overline{\otimes}
B(l^2))}\right\}\\
&=J\left(\sum_n dy_n^{(\nu)} \otimes e_{n,1}, 2^{\nu};
L^{1,\infty}(\M\overline{\otimes} B(l^2)),
L^{2}(\M\overline{\otimes} B(l^2))\right).
\end{split}
\end{equation*}
where $(e_{i,j})_{i,j}$ denotes the usual base of $B(l^2)$, that
is, $(dy_n^{(\nu)})_n$ is viewed as a column vector with entries
from $L^2(\M,\T)$.  This implies that
$$
\left\|\{ J(\sum_n dy_n^{(\nu)} \otimes e_{n,1},
2^{\nu})\}\right\|_{\lambda^{\theta,p}}\leq
K(\|x_\infty\|_{\theta,p;\underline{J}} +\epsilon).
$$
Let $S$ be a finite subset of $\mathbb{Z}$. By the definition of
$\|\cdot\|_{\theta,p;\underline{J}}$,
\begin{equation}\label{convergence1}
\left\|\sum_{\nu \in S} \sum_n dy_n^{(\nu)} \otimes e_{n,1}
\right\|_{[L^{1,\infty}(\M \overline{\otimes} B(l^2)) , L^{2}(\M
\overline{\otimes} B(l^2))]_{\theta,p;\underline{J}}} \leq
K(\|x_\infty\|_{\theta,p;\underline{J}} +\epsilon).
\end{equation}
Similar argument on $(z^{(\nu)})_\nu$ also gives,
\begin{equation}\label{convergence2}
\left\|\sum_{\nu \in S} \sum_n d{z_n^{(\nu)}}^* \otimes e_{n,1}
\right\|_{[L^{1,\infty}(\M \overline{\otimes} B(l^2)) , L^{2}(\M
\overline{\otimes} B(l^2))]_{\theta,p;\underline{J}}} \leq
K(\|x_\infty\|_{\theta,p;\underline{J}} +\epsilon).
\end{equation}
Set \[y:=\sum_{\nu\in\mathbb{Z}}y^{(\nu)}\ \ \text{and}\ \
z:=\sum_{\nu\in\mathbb{Z}}z^{(\nu)}.
\]
Note that since the  inequalities (\ref{convergence1}) and
(\ref{convergence2}) are  valid for arbitrary  finite subset $S$
of $\mathbb{Z}$, $y$ and $z$ are well-defined (that is, the series
are unconditionally convergent in the Banach space
$[L^{1,\infty}(\M,\T), L^{2}(\M,\T)]_{\theta,p;\underline{J}}$).
 Clearly, $y$ and $z$ are martingales   with $x=y +z$ and from
 (\ref{convergence1}) and
(\ref{convergence2}), we have:
\begin{equation}\label{BGinterpolation}
\left\|S_C(dy) \right\|_{\theta,p,\underline{J}} + \left\|S_R(dz)
\right\|_{\theta,p,\underline{J}}\leq
2K(\|x_\infty\|_{\theta,p;\underline{J}} +\epsilon).
\end{equation}
To conclude the proof, we note from the  general equivalence
theorem on real interpolations that the same statement  as in
(\ref{BGinterpolation}) can be made with any real interpolation
method (with possible change on the absolute constant). It is well
known that
\[
[L^{1,\infty}(\M,\T), L^2(\M,\T)]_{\theta,p}=L^p(\M,\T) \ \
\text{(with equivalent norms)}\]
 and
 \[
 [L^{1}(\M,\T),
L^2(\M,\T)]_{\theta,p}=L^p(\M,\T)\ \ \text{ (with equivalent
norms)}.
\]
 From \cite[Corollary~2.2, P.~1467]{PX3}, it is
enough to track the order of the constants for the commutative
case. Estimates on the order of growth of constants involved on
the equivalent norms can be deduced from \cite[Theorem~4.3]{Holm}
as follows: for $f\in L^2$,
\begin{equation}\label{commutative1}
C(1-\theta)^{-1/2}\|f\|_{L^p} \leq
\|f\|_{[L^{1,\infty},L^2]_{\theta, p}} \leq C^{-1}
\theta^{-1/p}(1-\theta)^{-1/p}\|f\|_{L^p}
\end{equation}
and
\begin{equation}\label{commutative2}
C\theta^{-1/p}(1-\theta)^{-1/2}\|f\|_{L^p} \leq
\|f\|_{[L^{1},L^2]_{\theta, p}} \leq C^{-1}
\theta^{-1}(1-\theta)^{-1/p}\|f\|_{L^p}.
\end{equation}

Combining (\ref{BGinterpolation}), (\ref{commutative1}), and
(\ref{commutative2}), we can conclude the existence of   an
absolute constant $M>0$ such that:
\[
\left\|S_C(dy) \right\|_p + \left\|S_R(dz) \right\|_p \leq
M(p-1)^{-1}(\|x\|_{p} +\epsilon).
\]
Taking the infimum over  $\epsilon>0$, we get
\[
 \|x\|_{\cal{H}^p(\M)} \leq M(p-1)^{-1}\|x\|_{p},
\]
which shows that $\alpha_p \leq M(p-1)^{-1}$ for $1<p<2$. The
proof is complete.
\end{proof}

\begin{remark} The optimal orders as stated in Theorem~\ref{alpha} remain valid
in the more general situation of Haagerup's $L^p$-spaces using
Haagerup's approximation (\cite{HAA2}). This follows from a
general deduction of the non-commutative Burkholder-Gundy
inequalities from the finite case to the type~III-case (with same
constants) achieved by Junge and Xu (still unpublished notes).
\end{remark}
For the next application, we recall the dual space of $\cal{H}^p$
($1<p<2$) studied in \cite{JX}. For $2<q\leq \infty$,
$L^q_CMO(\M)$ ($MO$ stands for {\em mean oscillation}) is the
space of all martingales  $x=(x_n)_n$ in $L^2(\M,\T)$ for which
\[
\|x\|_{L^q_CMO(\M)}=\sup_m\left\|\sup_{n\leq m}
\E_n\left(\sum_{k=n}^m |dx_k|^2\right)\right\|_{{q}/{2}}^{{1}/{2}}
<\infty.
\]
This was introduced as a non-commutative analogue of the $q$-norm
of the classical sharp function. In the above, the suggestive
notation introduced in \cite{Ju} for the supremum is understood in
the following sense: if $1\leq r, r'\leq \infty$ and $1/r +1/r'=1$
then for any sequence $(a_n)_{n\geq 1}$ of positive operators in
$\overline{\M}$,
\[
\left\|\sup_n a_n \right\|_r=\sup\left\{\sum_{n\geq 1}\T(a_nb_n):\
b_n\geq 0, \left\|\sum_{n\geq 1} b_n\right\|_{r'} \leq 1\right\}.
\]
Similarly, $L^q_RMO(\M)$ is defined as the space of all
martingales $x$ such that $x^*\in L^q_CMO(\M)$, with norm given by
$\|x\|_{L^q_RMO(\M)}=\|x^*\|_{L^q_CMO(\M)}$ and as above,
\[
L^qMO(\M)=L^q_CMO(\M) \cap L^q_RMO(\M),
\]
equipped with the usual
intersection norm
\[
\|x\|_{L^qMO(\M)}=\max\left\{\|x\|_{L^q_CMO(\M)},\|x\|_{L^q_RMO(\M)}\right\}.
\]
Let $1<p<2$ and $p'$ the index conjugate to $p$. It was shown in
\cite[Theorem~4.1]{JX} that $(\cal{H}^p(\M))^*=L^{p'}MO(\M)$. As
part of this characterization, they noted the  inequality that for
$2<q<\infty$, there is a constant $\lambda_q'>0$ such that
\[
\|a\|_q \leq \lambda_q'\|a\|_{L^qMO(\M)}, \quad a\in L^qMO(\M).
\]
 By duality, we can answer a problem raised in \cite[p.~972]{JX}  using Theorem~\ref{alpha}(i):
\begin{corollary}
$\lambda_q'\approx q$ as $q\to \infty$.
\end{corollary}

We end the paper with a short note on the class $L\log{L}$.
 Recall  the Zygmund space $L\log{L}$.
If $L^0(\Omega,\cal{F}, P)$ is the space of all (classes) of
measurable functions on a given probability space $(\Omega,
\cal{F}, P)$, the class $L\log{L}$ is defined by setting
$$L\log{L}=\left\{ f \in L^0(\Omega, \cal{F}, P);
\int |f|\log^+|f|\ dP <\infty\right\}.$$ Set $\|f\|_{L\log{L}}
=\int |f|\log^+|f|\ dP$. Note that $\|\cdot\|_{L\log{L}}$ is not a
norm but is  equivalent to a rearrangement invariant norm $||| f
|||_{L\log{L}} =\int_{0}^1 f^*(t) \log(1/t)\ dt$.

 Equipped with $|||\cdot|||_{L\log{L}}$,
 the spaces $L\log{L}$ is a rearrangement invariant Banach function
 space
(see for instance \cite[Theorem~6.4, pp.~246-247]{BENSHA}) so its
non-commutative analogue $L\log{L}(\M,\T)$   is well defined as
described in Section~1. We note as in \cite{Ran15} that if a
martingale $x$ is bounded in $L\log{L}(\M,\T)$ then it is
uniformly integrable in $L^1(\M,\T)$ and therefore is of the form
$x=(\E_n(x_\infty))_{n=1}^\infty$ for some  $x_\infty \in
L\log{L}(\M,\T)$.

 With $\alpha_p \approx (p-1)^{-1}$ (when $p\to 1$),  the elementary argument used in the proof of
\cite[Proposition~6.5]{Ran15} can be adjusted to   deduce the
following  strengthening of  \cite[Proposition~6.5]{Ran15}(which
answers positively \cite[Problem~6.4]{Ran15}).

\begin{theorem}\label{LlogL}
There is a constant $K>0$ such that if $x=(x_{n})^{\infty}_{n=1}$
is a martingale that is bounded in $L\log{L}(\M,\T)$, then
\begin{equation*}
\left\|x\right\|_{\cal{H}^1(\M)} \leq K + K\left\| x_\infty
\right\|_{L\log{L}(\M,\T)}.
\end{equation*}
\end{theorem}

\medskip

\noindent {\em Acknowledgements.} This project started during my
visit to the  Department of Mathematics at the Universit\'e
 de Franche-Comt\'e in Spring~2003. I would like to express my
gratitude to the department for  its  warm hospitality.



\providecommand{\bysame}{\leavevmode\hbox
to3em{\hrulefill}\thinspace}
\providecommand{\MR}{\relax\ifhmode\unskip\space\fi MR }
\providecommand{\MRhref}[2]{%
  \href{http://www.ams.org/mathscinet-getitem?mr=#1}{#2}
} \providecommand{\href}[2]{#2}

\end{document}